\documentclass[11pt]{article}
\usepackage{amsmath}
\usepackage{amssymb}
\usepackage{latexsym}
\usepackage{amsthm}
\def\disp{\displaystyle}

\def\F{Fr\'{e}chet}

\def\P{Painlev\'{e}}

\def\Limsup{\mathop{{\rm Lim}\,{\rm sup}}}

\def\tto{\;{\lower 1pt \hbox{$\rightarrow$}}\kern -10pt
\hbox{\raise 2pt \hbox{$\rightarrow$}}\;}
\def\lto{\longrightarrow}
\def\Hat{\widehat}
\def\Tilde{\widetilde}
\def\Bar{\overline}
\def\ra{\rangle}
\def\la{\langle}
\def\ve{\varepsilon}
\def\B{I\!\!B}
\def\h{\hfill\Box}\def\R{\mathbb{R}}
\def\N{\mathbb{N}}

\def\ox{\bar{x}}

\def\oz{\bar{z}}

\def\tx{\tilde{x}}

\def\epi{\mbox{\rm epi}\,}

\def\diam{\mbox{\rm diam}\,}

\def\cl{\mbox{\rm cl}\,}
\def\dist{\mbox{\rm dist}\,}

\def\cone{\mbox{\rm cone}\,}

\def\inter{\mbox{\rm int}\,}

\def\h{\hfill\square}
\def\dn{\downarrow}
\def\O{\Omega}

\def\ph{\varphi}
\def\emp{\emptyset}
\def\st{\stackrel}
\def\oR{\Bar{\R}}
\def\lm{\lambda}
\def\Lm{\Lambda}

\def\dd{\delta}
\def\al{\alpha}
\def\be{\beta}

\def\gg{\gamma}
\def \N{I\!\!N}

\def\vt{\vartheta}

\begin{document}
\begin{center}

{\bf RATED EXTREMAL PRINCIPLES FOR FINITE AND INFINITE
SYSTEMS}\footnote{This research was partially supported by the US
National Science Foundation under grants DMS-0603846 and
DMS-1007132 and by the Australian Research Council
under grant DP-12092508.}\\[3ex]
BORIS S. MORDUKHOVICH\footnote{Department of Mathematics, Wayne
State University, Detroit, MI 48202, USA. Email:
boris@math.wayne.edu.} and HUNG M. PHAN\footnote{Department of
Mathematics, Wayne State University, Detroit, MI 48202, USA.
Email: pmhung@wayne.edu.}\\[2ex]
{\bf Dedicated to Juan Enrique Martinez-Legaz in honor of his 60th
birthday}
\end{center}

\newtheorem{Theorem}{Theorem}[section]
\newtheorem{Proposition}[Theorem]{Proposition}
\newtheorem{Lemma}[Theorem]{Lemma}
\newtheorem{Corollary}[Theorem]{Corollary}
\newtheorem{Definition}[Theorem]{Definition}
\theoremstyle{definition}
\newtheorem{Remark}[Theorem]{Remark}
\newtheorem{Example}[Theorem]{Example}
\theoremstyle{plain}
\renewcommand{\theequation}{\thesection.\arabic{equation}}
\renewcommand{\theequation}{\thesection.\arabic{equation}}

\setcounter{equation}{0}\setcounter{section}{0}
\newcounter{count}
\newenvironment{lista}{\begin{list}{(\rm\alph{count})}{\usecounter{count}
\setlength{\itemindent}{-0.3cm} }}{\end{list}}

\newenvironment{listi}{\begin{list}{\rm(\textit{\roman{count}})}{\usecounter{count}
\setlength{\itemindent}{-0.3cm} }}{\end{list}}

\small{\bf Abstract.} In this paper we introduce new notions of
local extremality for finite and infinite systems of closed sets
and establish the corresponding extremal principles for them
called here rated extremal principles. These developments are in
the core geometric theory of variational analysis. We present
their applications to calculus and optimality conditions for
problems with infinitely many constraints. \vspace*{0.1cm}

\vspace*{0.1cm}


\section{Introduction}
\setcounter{equation}{0}

Modern variational analysis is based on variational principles and
techniques applied to optimization-related and equilibrium
problems as well as to a broad spectrum of problems, which may not
be of a variational nature; see the books \cite{Borwein-Zhu-TVA,
m-book1,m-book2,Rockafellar-Wets-VA} for more discussions and
references. In this vein, extremal principles have been well
recognize as fundamental geometric tools of variational analysis
and its applications that can be treated as far-going variational
extensions of convex separation theorems to systems of nonconvex
sets. We refer the reader to the two-volume monograph
\cite{m-book1,m-book2} and the bibliographies therein for various
developments and applications of the extremal principles in both
finite and infinite dimensions.

To the best of our knowledge, extremal principles have been
previously developed only for finite systems of sets. On the
other, there is a strong demand in various areas (e.g., in
semi-infinite optimization) for their counterparts involving
infinite, particularly countable, set systems.

The first attempt to deal with infinite systems of sets was
undertaken in our recent papers \cite{MorPh10a,MorPh10b}, where
certain tangential extremal principles were established for
countable set systems and then were applied therein to problems of
semi-infinite programming and multiobjective optimization. At the
same time, the tangential extremal principles developed and
applied in \cite{MorPh10a,MorPh10b} concern the so-called
tangential extremality (and only in finite dimensions) and do not
reduce to the conventional extremal principles of \cite{m-book1}
for finite systems of sets even in simple frameworks.

In this paper we develop new {\em rated extremal principles} for
both finite and infinite systems of closed sets in
finite-dimensional and infinite-dimensional spaces. Besides being
applied to conventional local extremal points of finite set
systems and reducing to the known results for them, the rated
extremal principles provide enhanced information in the case of
finitely many sets while open new lines of development for
countable set systems. The results obtained in this way allow us,
in particular, to derive intersection rules for generalized
normals of infinite intersections of closed sets, which imply in
turn new necessary optimality conditions for mathematical programs
with countable constraints in finite and infinite
dimensions.\vspace*{0.05in}

The rest of the paper is organized as follows. In Section~2 we
briefly discussed preliminaries from variational analysis and
generalized differentiations used in the sequel. In Section~3 we
introduce the notion of rated extremality and derive exact and
approximate versions of the rated extremal principles for systems
of finite sets in finite-dimensional and infinite-dimensional
spaces. Section~4 is devoted to rated extremal principles for
infinite/countable systems of closed sets in Banach spaces.
Finally, Section~5 provides applications of the rated extremal
principles to calculus of generalized normals to infinite set
intersections, which implies necessary optimality conditions for
optimization problems with countable geometric constraints.
\vspace*{0.05in}

Our notation is basically standard in variational analysis; see,
e.g., \cite{m-book1,Rockafellar-Wets-VA}. Recall that $B(\ox,r)$
stands for a closed ball centered at $\ox$ with radius $r>0$, that
$\B$ and $\B^*$ are the closed unit ball of the space in question
and its dual, respectively, and that $\N:=\{1,2,\ldots\}$. Given a
set-valued mapping $F\colon X\tto X^*$ between a Banach space $X$
and its topological dual $X^*$, we denote by
\begin{eqnarray}\label{eq:LimsupPK}
\begin{array}{ll}
\disp\Limsup_{x\to\ox} F(x):=\Big\{x^*\in X^*\in
Y\,\Big|&\exists\mbox{ sequences }\;x_k\to\ox\;\mbox{ and }\;
x^*_k\st{w^*}{\to}x^*\;\mbox{ as }\;k\to\infty\\
&\mbox{such that }\;x^*_k\in F(x_k)\mbox{ for all }\;k\in\N\Big\}
\end{array}
\end{eqnarray}
the {\em sequential Painlev\'{e}-Kuratowski outer limit} of $F$ at
$\ox$, where $w^*$ signifies the weak$^*$ topology of $X^*$.

\section{Preliminaries from Variational Analysis}
\setcounter{equation}{0}

In this section we briefly overview some basic tools of
variational analysis and generalized differentiation that are
widely used in what follows; see the books
\cite{Borwein-Zhu-TVA,m-book1,Rockafellar-Wets-VA,sc} for more
details and references. Unless otherwise stated, all the spaces
under consideration are Banach with the norm $\|\cdot\|$ and the
canonical pairing $\la\cdot,\cdot\ra$ between the space in
question and its topological dual.

Let $\O$ be a nonempty subset of a space $X$. Given $\ve\ge 0$,
the {\em set of $\ve$-normals} to $\O$ at $\ox$ is given by
\begin{equation}
\label{eq:e-nor}\Hat N_\ve(\ox;\O):=\left\{x^*\in
X^*\Big|\limsup_{x\st{\O}{\to}\ox}\frac{\la
x^*,x-\ox\ra}{\|x-\ox\|}\leq\ve\right\}
\end{equation}
with $\Hat N_\ve(\ox;\O):=\emp$ if $\ox\not\in\O$. When $\ve=0$,
the set \eqref{eq:e-nor} is denoted by $\Hat N(\ox;\O):=\Hat
N_0(\ox;\O)$ and is called the {\em \F\ normal cone} (or {\em
prenormal/regular normal cone}) to $\O$ at $\ox$. The {\em
Mordukhovich/basic/limiting normal cone} to $\O$ at a point
$\ox\in\O$ is defined by
\begin{equation}\label{eq:Lcone} N(\ox;\O):=\Limsup_{x\to\ox\atop\ve\dn0}\Hat
N_\ve(x;\O)
\end{equation}
via the sequential outer limit \P-Kuratowski outer limit
\eqref{eq:LimsupPK} of $\ve$-normals (\ref{eq:e-nor}) as $x\to\ox$
and $\ve\dn 0$. If the space $X$ is Asplund (i.e., each of its
separable subspace has a separable dual that holds, in particular,
when is reflexive) and the set $\O$ is locally closed around
$\ox$, we can equivalently put $\ve_k=0$ in \eqref{eq:Lcone}; see
\cite{m-book1} for more details. If $X=\R^n$, the basic normal
cone \eqref{eq:Lcone} can be equivalently described as
\begin{equation}\label{nc}
N(\ox;\O)=\Limsup_{x\to\ox}\Big\{\cone\big[x-\Pi(x;\O)\big]\Big\}
\end{equation}
via the Euclidian projector
$\Pi(x;\O):=\{w\in\O|\,\|x-w\|=\dist(x;\O)\}$ of $x\in\R^n$ onto
$\O$, which was the original definition in \cite{mor76}. In the
above formula \eqref{nc} the symbol $\cone A$ stands for the cone
generated by a nonempty set $A$ and is defined by
\begin{equation*}
\cone A:=\bigcup_{\lm\ge 0}\lm A.
\end{equation*}

Given an extended-real-valued function $\ph\colon
X\to\oR:=(-\infty,\infty]$, recall that the {\em Fr\'echet/regular
subdifferential} of $\ph$ at $\ox$ with $\ph(\ox)<\infty$ is
defined by
\begin{equation}\label{frechet}
\Hat\partial\ph(\ox):=\Big\{x^*\in
X^*\Big|\;\liminf_{x\to\ox}\frac{\ph(x)-\ph(\ox)-\la
x^*,x-\ox\ra}{\|x-\ox\|}\ge 0\Big\}.
\end{equation}
It is easy to see that $\Hat N(\ox;\O)=\Hat\partial\dd(\ox;\O)$
for the indicator function $\dd(\cdot;\O)$ of $\O$ defined by
$\dd(x;\O):=0$ when $x\in\O$ and $\dd(x;\O)=\infty$ otherwise.
Furthermore, we obviously have the following nonsmooth version of
the Fermat stationary rule:
\begin{equation}\label{fermat}
0\in\Hat\partial\ph(\ox)\;\mbox{ if }\;\ox\;\mbox{ is a local
minimizer of }\;\ph.
\end{equation}

A major motivation for our work is to develop and apply extremal
principles of variational analysis the first version of which was
formulated in \cite{KruMor80} for finitely many sets via
$\ve$-normals \eqref{eq:e-nor}; see \cite[Chapter~2]{m-book1} for
more details and discussions. Recall \cite[Definition
2.5]{m-book1} that a set system $\{\O_1,\ldots,\O_m\}$, $m\ge 2$,
satisfies the {\em approximate extremal principle} at
$\ox\in\cap_{i=1}^m\O_i$ if for every $\ve>0$ there are
$x_i\in\O_i\cap(\ox+\ve\B)$ and $x^*_i\in\Hat
N(x_i;\O_i)+\ve\B^*$, $i=1,\ldots,m$, such that
\begin{equation}
\label{eq:extr-n}
x^*_1+\ldots+x^*_m=0\quad\mbox{and}\quad\|x^*_1\|^2+\ldots+\|x^*_m\|^2=1.
\end{equation}
If the dual vectors $x^*_i$ can be taken from the limiting normal
cone $N(\ox;\O_i)$, then we say that the system
$\{\O_1,\ldots,\O_m\}$ satisfies the {\em exact extremal
principle} at $\ox$.

Efficient conditions ensuring the fulfillment of both approximate and exact versions of
the extremal principle can be found in \cite[Chapter~2]{m-book1} and the references
therein. Roughly speaking, the approximate extremal principle in terms of Fr\'echet
normals holds for locally extremal points of any closed subsets in Asplund spaces
(\cite[Theorem~2.20]{m-book1}) while the exact extremal principle requires additional
sequential normal compactness assumptions that are automatic in finite dimensions; see
\cite[Theorem~2.22]{m-book1}.

Recall \cite{KruMor80,m-book1} that a point
$\ox\in\cap_{i=1}^m\O_i$ is {\em locally extremal} for the system
$\{\O_1,\ldots,\O_m\}$ if there are sequences $\{a_{ik}\}\subset
X$, $i=1,\ldots,m$, and a neighborhood $U$ of $\ox$ such that
$a_{ik}\to 0$ as $k\to\infty$ and
\begin{equation}\label{eq:ESfin}
\bigcap_{i=1}^m\Big(\O_i-a_{ik}\Big)\cap U=\emp\;\mbox{ for all
large }\;k\in\N.
\end{equation}

As shown in \cite{m-book1}, this extremality notion for sets
encompasses standard notions of local optimality for various
optimization-related and equilibrium problems as well as for set
systems arising in proving calculus rules and other frameworks of
variational analysis.

\section{Rated Extremality of Finite Systems of Sets}
\label{sec:REPfin}\setcounter{equation}{0}

In this section we introduce a new notion of {\em rated
extremality} for finite systems of sets, which essentially broader
the previous notion \eqref{eq:ESfin} of local extremality. We show
nevertheless that both exact and approximate versions of the
extremal principle hold for this rated extremality under the same
assumptions as in \cite{m-book1} for locally extremal points. Let
us start with the definition of rated extremal points. For
simplicity we drop the word ``local" for rated extremal points in
what follows.

\begin{Definition}[Rated extremal points of finite set systems]
\label{Def:RES-fin} Let $\O_1,\ldots,\O_m$ as $m\ge 2$ be nonempty
subsets of $X$, and let $\ox$ be a common point of these sets. We
say that $\ox$ is a $($local$)$ {\sc rated extremal point} of rank
$\al$, $0\le\al<1$, of the set system $\{\O_1,\ldots,\O_m\}$ if
there are $\gg>0$ and sequences $\{a_{ik}\}\subset X$,
$i=1,\ldots,m$, such that $r_k:=\max_i \|a_{ik}\|\to 0$ as
$k\to\infty$ and
\begin{equation}\label{eq:RES-fin}
\bigcap_{i=1}^m\big(\O_i-a_{ik}\big)\cap B(\ox,\gg
r_k^\al)=\emptyset\quad\mbox{ for all large }\;k\in\N.
\end{equation}
In this case we say that $\{\O_1,\ldots,\O_m\}$ is a {\sc rated
extremal system} at $\ox$.
\end{Definition}

The case of local extremality \eqref{eq:ESfin} obviously
corresponds to \eqref{eq:RES-fin} with rate $\al=0$. The next
example shows that there are rated extremal points for systems of
two simple sets in $\R^2$, which are not locally extremal in the
conventional sense of \eqref{eq:ESfin}.

\begin{Example} [Rated extremality versus local extremality]
Consider the sets $$
\O_1:=\big\{(x_1,x_2)\in\R^2\big|\;x_2-x_1^2\le 0\big\}\;\mbox{
and }\;\O_2:=\big\{(x_1,x_2)\in\R^2\big|\;-x_2-x_1^2\le 0\big\}.
$$
Then it is easy to check that $(\ox_1,\ox_2)=(0,0)\in\O_1\cap\O_2$
is a rated extremal point of rank $\al=\frac{1}{2}$ for the system
$\{\O_1,\O_2\}$ but not a local extremal point of this system.
\end{Example}

Prior to proceeding with the main results of this section, we
briefly discuss relationships between the rated extremality and
the {\em tangential extremality} of set systems introduced in
\cite{MorPh10a}. Let $\{\O_i,i=1,\ldots,m\}$, $m\ge 2$, be  a
system of sets with $\ox\in\cap_{i=1}^m\O_i$, and let
$\Lm:=\{\Lm_i(\ox),i=1,\ldots,m\}$ be an approximating system of
cones. Recall that $\ox$ is a $\Lm$-tangential local extremal
point of $\{\O_i,i=1,\ldots,m\}$ if the system of cones
$\{\Lm_i(\ox),i=1,\ldots,m\}$ is extremal at the origin in the
sense that there are $a_1,\ldots,a_m\in X$ such that
\begin{equation*}
\bigcap_{i=1}^m\big(\O_i-a_i\big)=\emptyset.
\end{equation*}
We refer the reader to \cite{MorPh10a,MorPh10b} for more
discussion on the tangential extremality and its
applications.\vspace*{0.05in}

The next proposition result and the subsequent example reveal
relationships between the rated extremality and tangential
extremality of set systems.

\begin{Proposition} [Relationships between rated and tangential
extremality of finite systems of sets] \label{Thm:Tan-Rated} Let
$\{\O_1,\ldots,\O_m\}$ as $m\ge 2$ be a $\Lm$-tangential extremal
system of sets at $\ox$. Assume that there are real numbers $C>0$,
$p\in(0,1)$ and a neighborhood $U$ of $\ox$ such that
\begin{equation}\label{t1}
\dist(x-\ox;\Lm_i)\le C\|x-\ox\|^{1+p}\;\mbox{ for all }\;
x\in\O_i\cap U\;\mbox{ and }\;i=1,\ldots,m.
\end{equation}
Then $\{\O_1,\ldots,\O_m\}$ is a rated extremal system at $\ox$.
\end{Proposition}
{\bf Proof.} Since the general case of $m\ge 2$ can be derived by
induction, it suffices to justify the result in the case of $m=2$.
Let $\{\Lm_1,\Lm_2\}$ be an extremal system of approximation cones
and find by definition elements $a_1,a_2\in X$ such that
\begin{equation*}
(\Lm_1-a_1)\cap(\Lm_2-a_2)=\emp.
\end{equation*}
Without loss of generality, assume that $a_1=-a_2=:a$. Take
$\al\in(0,1)$ with $\be:=\al(1+p)>1$ and show that for all small
$t>0$ we have
\begin{equation}
(\O_1-ta)\cap(\O_2+ta)\cap B(\ox,\|ta\|^\al)=\emp.
\end{equation}
Suppose by contradiction that there exists
\begin{equation}
x\in(\O_1-ta)\cap(\O_2+ta)\cap B(\ox,\|ta\|^\al).
\end{equation}
That implies by using condition \eqref{t1} that
\begin{equation*}
\begin{aligned}
&\dist(x-\ox;\Lm_1-ta)=\dist(x+ta-\ox;\Lm_1)\le C\|x+ta-\ox\|^{1+p},\\
&\dist(x-\ox;\Lm_2+ta)=\dist(x-ta-\ox;\Lm_2)\le
C\|x-ta-\ox\|^{1+p}.
\end{aligned}
\end{equation*}
Thus we have for some constant $\Tilde C$ that
\begin{equation*}
\|x+ta-\ox\|^{1+p}\le\Tilde
C\max\big\{\|x-\ox\|,\|ta\|\big\}^{1+p}\le \Tilde
C\max\big\{\|ta\|^\be,\|ta\|^{1+p}\big\}=o(\|ta\|)\;\mbox{ as }\;
t\dn 0
\end{equation*}
and similarly $\|x-ta-\ox\|^{1+p}=o(\|ta\|)$. Put then
$d:=\dist(\Lm_1-a,\Lm_2+a)>0$ and observe due the conic structures
of $\Lm_1$ and $\Lm_2$ that
\begin{equation*}
td=\dist(\Lm_1-ta;\Lm_2+ta)>0
\end{equation*}
for all $t>0$ sufficiently small. Combining all the above gives us
\begin{equation*}
td=\dist(\Lm_1-ta;\Lm_2+ta)\le
\dist(x-\ox;\Lm_1-ta)+\dist(x-\ox;\Lm_2+ta)=o(\|ta\|),
\end{equation*}
which is a contradiction. Thus $\{\O_1,\O_2,\ox\}$ is a rated
extremal system at $\ox$ with rank $\al$ chosen above. This
completes the proof of the proposition. $\h$\vspace*{0.05in}

One of the most important special cases of tangential extremality
is the so-called {\em contingent extremality} when the
approximating cones to $\O_i$ are given by the Bouligand-Severi
contingent cones to this sets; see \cite{MorPh10a,MorPh10b}, where
this case of tangential extremality was primarily studied and
applied. The following example (of two parts) shows that the
notions of rated extremality and contingent extremality are
independent from each other in a simple setting of two sets in
$\R^2$.

\begin{Example} [Independence of rated and contingent
extremality] Let $X=\R^2$, and let $\ox=(0,0)$.

{\bf (i)} Consider two closed sets in $\R^2$ given by
\begin{equation*}
\O_1:=\epi f\;\mbox{ and }\;\O_2:=\R\times\R_-\setminus\inter\O_1,
\end{equation*}
where $f(x):=x\sin\frac{1}{x}$ for $x\in\R$ with $f(0):=0$. It is
easy to see that the contingent cones to $\O_1$ and $\O_2$ at
$\ox$ are computed by
\begin{equation*}
\Lm_1=\epi(-|\cdot|)\;\mbox{ and }\;\Lm_2=\R\times\R_-.
\end{equation*}
We can check that the set system $\{\O_1,\O_2\}$ is locally
extremal at $\ox$, and hence $\ox$ is a rated extremal point of
this system of sets with rank $\al=0$. On the other hand, the
contingent extremality is obviously violated for $\{\O_1,\O_2\}$
at $\ox$ as follows from the above computations of $\Lm_1$ and
$\Lm_2$.

{\bf (ii)} Now we define two closed sets in $\R^2$ by
\begin{equation*}
\O_1:=\R\times\R_-\;\mbox{ and }\;\O_2:=\epi\;\mbox{ with
}\;f(x):=-x^{1+\tfrac{1}{\ln^2|x|}}\;\mbox{ for }\;x\ne 0\;\mbox{
and }\;f(0):=0.
\end{equation*}
The contingent cones to $\O_1$ and $\O_2$ at $\ox$ are easily
computed by $\Lm_1=\R\times\R_-$ and  $\Lm_2=\R\times\R_+$. We can
check that $\ox$ is not a rated extremal point of $\{\O_1,\O_2\}$
whenever $\al\in[0,1)$, while the contingent extremality obviously
holds for this system at $\ox$.
\end{Example}

The next theorem justifies the fulfillment of the exact extremal
principle for any rated extremal point of a finite system of
closed sets in $\R^n$. It extends the extremal principle of
\cite[Theorem~2.8]{m-book1} obtained for local extremal points,
i.e., when $\al=0$ in Definition~\ref{Def:RES-fin}.

\begin{Theorem}[Exact extremal principle for rated extremal systems of
sets in finite dimensions] Let $\ox$ be a rated extremal point of
rank $\al\in[0,1)$ for the system of sets $\{\O_1,\ldots,\O_m\}$
as $m\ge 2$ in $\R^n$. Assume that all the sets $\O_i$ are locally
closed around $\ox$. Then the exact extremal principle holds for
$\{\O_1,\ldots,\O_m\}$ at $\ox$, i.e, there are $x^*_i\in
N(\ox;\O_i)$ for $i=1,\ldots,m$ satisfying the relationships in
\eqref{eq:extr-n}.
\end{Theorem}
{\bf Proof.} Given a rated extremal point $\ox$ of the system
$\{\O_1,\ldots,\O_m\}$, take numbers $\al\in[0,1)$ and $\gg>0$ as
well as sequences $\{a_{ik}\}$ and $\{r_k\}$ from
Definition~\ref{Def:RES-fin}. Consider the following unconstrained
minimization problem for any fixed $k\in\N$:
\begin{equation}
\label{eq:P1k} \mbox{minimize }\;
d_k(x):=\left[\sum_{i=1}^m\dist^2\big(x+a_{ik};\O_i\big)\right]^\frac{1}{2}+\frac{\sqrt{m}}
{\gg^{\frac{1}{\al}}}\|x-\ox\|^{\frac{1}{\al}}\ ,\ x\in\R^n.
\end{equation}
Since the function $d_k$ is continuous and its level sets are
bounded, there exists an optimal solution $x_k$ to (\ref{eq:P1k})
by the classical Weierstrass theorem. We obviously have the
relationships
\begin{equation*}
d_k(x_k)\le
d_k(\ox)=\left[\sum_{i=1}^m\dist^2\big(\ox+a_{ik};\O_i\big)\right]^\frac{1}{2}\le
\left[\sum_{i=1}^m\|a_{ik}\|^2\right]^\frac{1}{2}\le r_k\sqrt{m},
\end{equation*}
which readily imply the estimate
\begin{equation*}
\displaystyle\frac{\sqrt{m}}{\gg^{\frac{1}{\al}}}\|x_k-\ox\|^{\frac{1}{\al}}\le
r_k\sqrt{m},\;\mbox{ i.e., }\;\|x_k-\ox\|\le\gg r_k^{\al}.
\end{equation*}
Taking the latter into account, we get
\begin{equation*}
\nu_k:=\left[\sum_{i=1}^m\dist^2\big(x_k+a_{ik};\O_i\big)\right]^\frac{1}{2}>0,
\end{equation*}
since the opposite statement $\nu_k=0$ contradicts the rated
extremality of $\ox$. Furthermore, the optimality of $x_k$ in
\eqref{eq:P1k} and choice of $\{a_{ik}\}$ give us the
relationships
\begin{equation*}
d_k(x_k)=\nu_k+\frac{\sqrt{m}}{\gg^{\frac{1}{\al}}}\|x_k-\ox\|^{\frac{1}{\al}}\le\left[\sum_{i=1}^m
\|a_{ik}\|^2\right]^\frac{1}{2}\downarrow 0\;\mbox{ as
}\;k\to\infty,
\end{equation*}
which ensure in turn that $x_k\to\ox$ and $\nu_k\downarrow 0$ as
$k\to\infty$.

We now arbitrarily pick $w_{ik}\in\Pi(x_k+a_{ik};\O_i)$ for
$i=1,\ldots,m$ in the closed set $\O_i$ and for each $k\in\N$
consider the problem:
\begin{equation}
\label{eq:P2k} \mbox{minimize }\;
\rho_k(x):=\left[\sum_{i=1}^m\|x+a_{ik}-w_{ik}\|^2\right]^\frac{1}{2}
+\frac{\sqrt{m}}{\gg^{\frac{1}{\al}}}\|x-\ox\|^{\frac{1}{\al}},\quad
x\in\R^n,
\end{equation}
which obviously has the same optimal solution $x_k$ as for
(\ref{eq:P1k}). Since $\nu_k>0$ and the norm $\|\cdot\|$ is
Euclidian, the function $\rho_k(\cdot)$ in \eqref{eq:P2k} is
continuously differentiable around $x_k$. Thus applying the
classical Fermat rule to the {\em smooth} unconstrained
minimization problem (\ref{eq:P2k}), we get
\begin{equation*}
\nabla\rho_k(x_k)=\sum_{i=1}^m
x^*_{ik}+C\|x_k-\ox\|^{\frac{1-2\al}{\al}}(x_k-\ox)=0\quad{\rm
for\ some\ constant\ }C,
\end{equation*}
where $x^*_{ik}:=(x_k+a_{ik}-w_{ik})/\nu_k$ for $i=1,\ldots,m$
with
\begin{equation*}
\|x^*_{1k}\|^2+\ldots+\|x^*_{mk}\|^2=1.
\end{equation*}
Observe that
$\displaystyle\|x_k-\ox\|^{\frac{1-2\al}{\al}}(x_k-\ox)=
\|x_k-\ox\|^{\frac{1-\al} {\al}}\frac{x_k-\ox}{\|x_k-\ox\|}\to 0$
as $x_k\to \ox$. Due to the compactness of the unit sphere in
$\R^n$, we find $x^*_i\in\R^n$ as $i=1,\ldots,m$ such that
$x^*_{ik}\to x^*_i$ as $k\to\infty$ without relabeling. It follows
from the equivalent description \eqref{nc} of the limiting normal
cone that $x^*_i\in N(\ox;\O_i)$ for all $i=1,\ldots,m$. Moreover,
we get from the constructions above that
\begin{align*}
\|x^*_{1}\|^2+\ldots+\|x^*_{m}\|^2=1\;\mbox{ and }\;
x^*_1+\ldots+x^*_m=0.
\end{align*}
This gives all the conclusions of the exact extremal principle and
completes the proof of the theorem. $\h$\vspace*{0.05in}

The next example shows that the exact extremal principle is
violated if we take $\al=1$ in Definition~\ref{Def:RES-fin}.

\begin{Example}[Violating the exact extremal principle for rated extremal points
of rank $\al=1$]\label{viol} Define two closed sets in $\R^2$ by
\begin{equation*}
\O_1:=\epi(-\|\cdot\|)\;\mbox{ and }\;\O_2:=\R\times\R_-.
\end{equation*}
Taking any $a_k\downarrow 0$, we see that
\begin{equation*}
\big(\O_1+(0,a_k)\big)\cap\big(\O_1-(0,a_k)\big)\cap
B(\ox,a_k/2)=\emptyset,
\end{equation*}
i.e., $\ox=(0,0)$ is a rated extremal point of $\{\O_1,\O_2\}$ of
rank $\al=1$. However, it is easy to check that the relationships
of the exact extremal principle do not hold for this system at
$\ox$.
\end{Example}

Observe that Example~\ref{viol} shows that the relationships of
the approximate extremal principle are also violated when $\al=1$.
However, for rated extremal systems of rank $\al\in[0,1)$ the
approximate extremal principle holds in general
infinite-dimensional settings. Let us proceed with justifying this
statement extending the corresponding results of \cite{m-book1}
obtained for the rank $\al=0$ in Definition~\ref{Def:RES-fin}.

\begin{Theorem}[Approximate extremal principle for rated extremal systems in
\F\ smooth spaces] Let $X$ be a Banach space admitting an
equivalent norm Fr\'echet differentiable off the origin, and let
$\ox$ be a rated extremal point of rank $\al\in[0,1)$ for a system
of sets $\O_1,\ldots,\O_m$ locally closed around $\ox$. Then the
approximate extremal principle holds for $\{\O_1,\ldots,\O_m\}$ at
$\ox$.
\end{Theorem}
{\bf Proof.} Choose an equivalent norm $\|\cdot\|$ on $X$
differentiable off the origin and consider first the case of $m=2$
in the theorem. Let $\ox\in\O_1\cap\O_2$ be a rated extremal point
of rank $\al\in[0,1)$ with $\gg>0$ taken from
Definition~\ref{Def:RES-fin}. Denote $r:=\max\{\|a_1\|,\|a_2\|\}$
and for any $\ve>0$ find $a_1,a_2$ such that
\begin{equation*}
r^{1-\al}\le\min\Big\{\frac{\gg}{2},\frac{\ve}{(2\gg
)^{(1-\al)/\al}}\Big\}\;\mbox{ and
}\;\big(\O_1-a_1\big)\cap\big(\O_2-a_2\big)\cap B\big(\ox,\gg
r^{\al}\big)=\emptyset.
\end{equation*}
We also select a constant $C>0$ with
$(\frac{2}{C})^{\al}=\frac{\gg}{2}$ and denote
$\be:=\frac{1}{\al}>1$. Define the function
\begin{equation}\label{ph}
\ph(z):=\|(x_1-a_1)-(x_2-a_2)\|\quad\mbox{for}\quad z=(x_1,x_2)\in
X\times X
\end{equation}
with the product norm $\|z\|:=(\|x_1\|^2+\|x_2\|^2)^{1/2}$ on
$X\times X$, which is \F\ differentiable off the origin under this
property of the norm on $X$. Next fix $z_0=(\ox,\ox)$ and define
the set
\begin{equation}\label{w0}
W(z_0):=\left\{z\in\O_1\times\O_2\big|\ph(z)+C\|z-z_0\|^\be\leq\ph(z_0)\right\},
\end{equation}
which is obviously nonempty and closed. For each $z=(x_1,x_2)\in
W(z_0)$ we have $i=1,2$:
\begin{equation*}
C\|x_i-\ox\|^\be\le C\|z-\oz\|^\be\le\ph(z_0)=\|-a_1+a_2\|\le
2r,\quad i=1,2,
\end{equation*}
which implies that $\displaystyle \|x_i-\ox\|\le
\left(\tfrac{2}{C}\right)^\frac{1}{\be}r^\frac{1}{\be}=
\left(\tfrac{2}{C}\right)^{\al}r^{\al}=\tfrac{\gg }{2}r^{\al}$ and
thus
\begin{equation*}
W(z_0)\subset B(\ox,\gg r^{\al})\times B(\ox,\gg r^{\al})\subset
B\big(\ox,\tfrac{1}{2}\ve^{\frac{\al}{1-\al}}\big)\times
B\big(\ox,\tfrac{1}{2}\ve^{\frac{\al}{1-\al}}\big).
\end{equation*}
It follows from Definition~\ref{Def:RES-fin} and constructions
\eqref{ph} and \eqref{w0} that $\ph(z)>0$ for all $z\in W(x_0)$.
Indeed, assuming on the contrary that $\ph(z)=0$ for some
$z=(x_1,x_2)\in W(x_0)$ gives us
\begin{equation*}
\|x_1-a_1-\ox\|\le\|x_1-\ox\|+\|a_1\|\le\tfrac{\gg }{2}r^\al+r=
\left(\tfrac{\gg}{2}+r^{1-\al}\right)r^\al\le\gg r^\al
\end{equation*}
and thus
$x_1-a_1=x_2-a_2\in\big(\O_1-a_1\big)\cap\big(\O_2-a_2\big)\cap
B(\ox,\gg r^\al)\ne\emp$, a contradiction.\vspace*{0.05in}

Hence $\ph$ is \F\ differentiable at any point $z\in W(z_0)$. Pick
any $z_1\in\O_1\times\O_2$ satisfying
\begin{equation*}
\ph(z_1)+C\|z_1-z_0\|^\be\le\inf_{W(z_0)}\Big\{\ph(z)+C\|z-z_0\|^\be\Big\}+\frac{r}{2}
\end{equation*}
and define further the nonempty and closed set
\begin{equation*}
W(z_1):=\left\{z\in\O_1\times\O_2\Big|\;\ph(z)+C\|z-z_0\|^\be+C\frac{\|z-z_1\|^\be}{2}
\le\ph(z_1)+C\|z_1-z_0\|^\be\right\}.
\end{equation*}
Arguing inductively, suppose we have chosen $z_k$ and constructed
$W(z_k)$, then pick $z_{k+1}\in W(z_k)$ such that
\begin{equation*}
\ph(z_{k+1})+C\sum_{i=0}^k\frac{\|z_{k+1}-z_i\|^\be}{2^i}\le
\inf_{W(z_k)}\Big\{\ph(z)+C\sum_{i=0}^k\frac{\|z-z_i\|^\be}{2^i}\Big\}+\frac{r}{2^{2k+1}}
\end{equation*}
and construct the subsequent nonempty and closed set
\begin{equation*}
W(z_{k+1}):=\left\{z\in\O_1\times\O_2\Big|\;\ph(z)+C\sum_{i=0}^{k+1}\frac{\|z-z_i\|^\be}{2^i}
\le\ph(z_{k+1})+C\sum_{i=0}^{k}\frac{\|z_{k+1}-z_i\|^\be}{2^i}
\right\}.
\end{equation*}
It is easy to see that the sequence
$\{W(z_k)\}\subset\O_1\times\O_2$ is nested. Let us check that
\begin{equation}\label{diam}
\diam W(z_{k+1}):=\sup\big\{\|z-w\|\,\big|\,z,w\in
W(z_{k+1})\big\}\to 0\;\mbox{  as }\;k\to\infty.
\end{equation}
Indeed, for each $z\in W(z_{k+1})$ and $k\in\N$ we have
\begin{align*}
C\frac{\|z-z_{k+1}\|^\be}{2^{k+1}}&\le
\ph(z_{k+1})+C\sum_{i=0}^{k}\frac{\|z_{k+1}-z_i\|^\be}{2^i}-
\left(\ph(z)+C\sum_{i=0}^{k}\frac{\|z-z_i\|^\be}{2^i}\right)\\
&\le\ph(z_{k+1})+C\sum_{i=0}^{k}\frac{\|z_{k+1}-z_i\|^\be}{2^i}-
\inf_{W(z_k)}\Big\{\ph(z)+C\sum_{i=0}^k\frac{\|z-z_i\|^\be}{2^i}\Big\}
\le\frac{r}{2^{2k+1}},
\end{align*}
which implies that $\displaystyle\diam W(z_{k+1})\le
2\left(\frac{r}{C2^k}\right)^\frac{1}{\be}$ and thus justifies
\eqref{diam}. Due to the completeness of $X$ the classical Cantor
theorem ensures the existence of $\oz=(\ox_1,\ox_2)\in W(z_0)$
such that $\displaystyle\bigcap_{k=0}^{\infty}W(z_k)=\{\oz\}$ with
$z_k\to\oz$ as $k\to\infty$. Now we show that $\oz$ is a minimum
point of the function
\begin{equation}\label{phi}
\phi(z):=\ph(z)+C\sum_{i=0}^\infty\frac{\|z-z_i\|^\be}{2^i}
\end{equation}
over the set $\O_1\times\O_2$. To proceed, take any $\oz\ne
z\in\O_1\times\O_2$ and observe that $z\not\in W(z_k)$ for all
$k\in\N$ sufficiently large while $\oz\in W(z_k)$. This yields the
estimates
\begin{equation*}
\phi(z)\ge\ph(z)+C\sum_{i=0}^k\frac{\|z-z_i\|^\be}{2^i}\ge
\ph(z_k)+C\sum_{i=0}^{k-1}\frac{\|z_k-z_i\|^\be}{2^i}\ge
\ph(\oz)+C\sum_{i=0}^k\frac{\|\oz-z_i\|^\be}{2^i}
\end{equation*}
and hence justifies the claimed inequality $\phi(z)\ge\phi(\oz)$
by letting $k\to\infty$.

We get therefore that the function
$\phi(z)+\delta(z;\O_1\times\O_2)$ attains at $\oz$ its minimum on
the whole space $X\times X$. The generalized Fermat rule
\eqref{fermat} gives us the inclusion
$0\in\Hat\partial\big(\phi(z)+\delta(z;\O_1\times\O_2)\big)$.
Since $\ph(\oz)>0$ and the norm $\|\cdot\|^\be$ is smooth, the
function $\phi$ in \eqref{phi} is \F\ differentiable at $\oz$.
Applying the sum rule from \cite[Proposition~1.107]{m-book1}, the
\F\ subdifferential formula for the indicator function, and the
product formula for Fr\'echet normal cone \eqref{eq:e-nor} from
\cite[Proposition~1.2]{m-book1}, we get
\begin{equation*}
-\nabla\phi(\oz)=-(u^*_1,u^*_2)\in\Hat N(\oz;\O_1\times\O_2)=\Hat
N(\ox_1;\O_1)\times\Hat N(\ox_2;\O_2),
\end{equation*}
where the dual elements $u^*_i$, $i=1,2$, are computed by
\begin{equation*}
u^*_1=x^*+\sum_{j=0}^\infty
w^*_{1j}\frac{\|\ox_1-x_{1j}\|^{\be-1}}{2^j}\;\mbox{ and }\;
u^*_2=-x^*+\sum_{j=0}^\infty
w^*_{2j}\frac{\|\ox_2-x_{2j}\|^{\be-1}}{2^j}
\end{equation*}
with $z_j=(x_{1j},x_{2j})$,
$x^*=\nabla\big(\|\cdot\|\big)\big((\ox_1-a_1)-(\ox_2-a_2)\big)$,
and
\begin{equation*}
w^*_{ij}=
\begin{cases}
\nabla(\|\cdot\|)(\ox_i-x_{ij})&\mbox{if }\;\ox_i-x_{ij}\ne 0,\\
0&\mbox{otherwise}.
\end{cases}
\end{equation*}
for $i=1,2$ and $j=0,1,\ldots$ due to the construction of the
function $\phi$ in \eqref{phi}. Observing further that $\|x^*\|=1$
and that $\oz,z_i\in W(z_0)$ gives us
\begin{equation*}
\|\ox_i-x_{ij}\|\le\ve^\frac{1-\al}{\al}=\ve^\frac{1}{\be-1},
\end{equation*}
which implies the estimates $\|\ox_i-x_{ij}\|^{\be-1}\le\ve$ and
\begin{equation*}
\sum_{j=0}^\infty\|w^*_{ij}\|\frac{\|\ox_i-x_{ij}\|^{\be-1}}{2^j}\le
2\ve,\quad i=1,2.
\end{equation*}
Setting finally $x^*_1:=-x^*/2$, $x^*_2:=x^*/2$, and $x_i:=\ox_i$
for $i=1,2$, we arrive at the relationships
\begin{align*}
&x^*_i\in\Hat N(x_i;\O_i)+\ve B^*\ ,\ x_i\in B(\ox,\ve)\;\mbox{
for }\;i=1,2,\\ &\|x^*_1\|+\|x^*_2\|=1,\;\mbox{ and
}\;x^*_1+x^*_2=0,
\end{align*}
which show that the approximate extremal principle holds for rated
extremal points of two sets.

Consider now the general case of $m>2$ sets. Observe that if $\ox$
as a rated extremal point of the system $\{\O_1,\ldots,\O_m\}$
with some rank $\al\in[0,1)$, then the point
$\oz:=(\ox,\ldots,\ox)\in X^{n-1}$ is a local rated extremal point
of the same rank for the system of two sets
\begin{equation}\label{m1}
\Theta_1:=\O_1\times\ldots\times\O_{n-1}\ \mbox{ and }\
\Theta_2:=\big\{(x,\ldots,x)\in X^{n-1}\big|x\in\O_m\big\}.
\end{equation}
To justify this, take numbers $\al\in[0,1)$ and $\gg >0$ and the
sequences $(a_{1k},\ldots,a_{mk})$ from
Definition~\ref{Def:RES-fin} for $m$ sets and check that
\begin{equation}\label{m2}
\Big(\Theta_1-(a_{1k},\ldots,a_{n-1,k})\Big)\cap\Big(\Theta_2-(a_{nk},\ldots,a_{nk})
\Big)\cap B\big((\ox,\ldots,\ox);\gg r_k^\al\big)=\emp
\end{equation}
with $r_k:=\max\{\|a_{1k}\|,\ldots,\|a_{nk}\|\}$. Indeed, the
violation of \eqref{m2} means that there are
$(x_1,\ldots,x_{n-1})\in\O_1\times\ldots\times\O_{n-1}$ and
$x_m\in\O_m$ satisfying
\begin{equation*}
x_1-a_{1k}=\ldots=x_{m-1}-a_{m-1,k}=x_m-a_{mk}\in B(\ox,\gg
r_k^\al),
\end{equation*}
which clearly contradicts the rated extremality of $\ox$ with rank
$\al$ for the system $\{\O_1,\ldots,\O_m\}$. Applying finally the
relationships of the approximate extremal principle to the system
of two sets in \eqref{m1} and taking into account the structures
of these sets as well as the aforementioned product formula for
Fr\'echet normals, we complete the proof of the theorem.
$\h$\vspace*{0.05in}

The next theorem elevates the fulfillment of the approximate
extremal principle for rated extremal points from Fr\'echet smooth
to Asplund spaces by using the method of {\em separable
reduction}; see \cite{Fab-Mor02,m-book1}.

\begin{Theorem}[Approximate extremal principle for rated extremal systems in Asplund
spaces]\label{aep} Let $X$ be an Asplund space, and let $\ox$ be a
rated extremal point of rank $\al\in[0,1)$ for a system of sets
$\O_1,\ldots,\O_m$ locally closed around $\ox$. Then the
approximate extremal principle holds for $\{\O_1,\ldots,\O_m\}$ at
$\ox$.
\end{Theorem}
{\bf Proof.} Taking a rated extremal point $\ox$ for the system
$\{\O_1,\ldots,\O_m\}$ of rank $\al\in[0,1)$, find a number
$\gg>0$ and sequences $\{a_{ik}\}$, $i=1,\ldots,m$, from
Definition~\ref{Def:RES-fin}. Consider a separable subspace $Y_0$
of the Asplund space $X$ defined by
\begin{equation*}
Y_0:=\mbox{span}\big\{\ox,a_{ik}\big|\;i=1,\ldots,m,\
k\in\N\big\}.
\end{equation*}
Pick now a closed and separable subspace $Y\subset X$ with
$Y\supset Y_0$ and observe that $\ox$ is a rated extremal point of
rank $\al$ for the system $\{\O_1\cap Y,\ldots,\O_m\cap Y\}$.
Indeed, we have
\begin{align*}
&\Big((\O_1\cap Y)-a_{1k}\Big)\cap\ldots\cap\Big((\O_m\cap Y)-a_{mk}\Big)
\cap B_Y(\ox;\gg r_k^\al)\\
&\subset
\Big(\O_1-a_{1k}\Big)\cap\ldots\cap\Big(\O_m-a_{mk}\Big)\cap
B_X(\ox;\gg r_k^\al)=\emptyset,
\end{align*}
where $r_k:=\max\{\|a_{1k}\|,\ldots,\|a_{mk}\|\}$, and where $B_X$
and $B_Y$ are the closed unit balls in the space $X$ and $Y$,
respectively. The rest of the proof follows the one in
\cite[Theorem 2.20]{m-book1} by taking into account that $Y$
admits an equivalent Fr\'echet differentiable norm off the origin.
$\h$ \vspace{0.05in}

We conclude this section with deriving the exact extremal
principle for rated extremal systems of rank $\al\in[0,1)$ in
Asplund spaces extending the corresponding result of
\cite[Theorem~2.22]{m-book1} obtained for $\al=0$.

Recall that a set $\O\subset X$ is {\em sequentially normally
compact} (SNC) at $\ox\in\O$ if for any sequence
$\{(x_k,x^*_k)\}_{k\in\N}\subset\O\times X^*$ we have the
implication
\begin{equation}\label{snc}
\big[x_k\to\ox,\;x^*_k\st{w^*}{\to}0\;\mbox{ with }\;x^*_k\in\Hat
N(x_k;\O),\;k\in\N\big]\Longrightarrow\|x^*_k\|\to
0\; \mbox{ as }\;k\to\infty.
\end{equation}
Besides the obvious validity of this property in
finite-dimensional spaces, it holds also in broad
infinite-dimensional settings; see, in particular,
\cite[Subsection~1.2.5]{m-book1} and SNC calculus rules
established in \cite[Section~3.3]{m-book1} in the framework of
Asplund spaces.

\begin{Theorem}[Exact extremal principle for rated extremal systems in Asplund
spaces]\label{eep} Let $X$ be an Asplund space, and let $\ox$ be a
rated extremal point of rank $\al\in[0,1)$ for a system of sets
$\O_1,\ldots,\O_m$ locally closed around $\ox$. Assume that all
but one of the sets $\O_i$, $i=1,\ldots,m$, are SNC at $\ox$. Then
the exact extremal principle holds for $\{\O_1,\ldots,\O_m\}$ at
$\ox$.
\end{Theorem}
{\bf Proof.} Follows the lines in the proof of
\cite[Theorem~2.22]{m-book1} by passing to the limit in the
relationships of the rated approximate extremal principle obtained
in Theorem~\ref{aep}. $\h$

\section{Rated Extremal Principles for Infinite Set Systems}
\label{sec:REPinf} \setcounter{equation}{0}

This section concerns new notions of rated extremality and
deriving rated extremal principles for infinite systems of closed
sets. The main results are obtained in the framework of Asplund
spaces.\vspace*{0.05in}

Let us start with introducing a notion of rated extremality for
arbitrary (may be infinite and not even countable) systems of sets
in general Banach spaces. We say that $R(\cdot)\colon\R_+\to\R_+$
is a {\em rate function} if there is a real number $M$ such that
\begin{equation}\label{rate}
rR(r)\le M\;\mbox{ and }\;\lim_{r\dn 0}R(r)=\infty.
\end{equation}
In what follow we denote by $|I|$ the cardinality (number of
elements) of a finite set $I$.

\begin{Definition}[Rated extremality for infinite systems of sets]
\label{Def:RES-Inf} Let $\{\O_i\}_{i\in T}$ be a system of closed
subsets of $X$ indexed by an arbitrary set $T$, and let
$\ox\in\bigcap_{t\in T}\O_i$. Given a rate function $R(\cdot)$, we
say that $\ox$ is an $R$-{\sc rated extremal point} of the system
$\{\O_i\}_{i\in T}$ if there exist sequences $\{a_{ik}\}\subset
X$, $i\in T$ and $k\in\N$, with $r_k:=\sup_{i\in T}\|a_{ik}\|\to
0$ as $k\to\infty$ such that whenever $k\in\N$ there is a finite
index subset $I_k\subset T$ of cardinality $|I_k|^{3/2}=o(R_k)$
with $R_k:=R(r_k)$ satisfying
\begin{equation}\label{eq:RES-Inf}
\bigcap_{i\in I_k}\big(\O_i-a_{ik}\big)\cap
B\big(\ox;r_kR_k\big)=\emp\;\mbox{ for all large }\;k.
\end{equation}
In this case we say that $\{\O_i\}_{i\in T}$ is an $R$-{\sc rated
extremal system} at $\ox$.
\end{Definition}

It is easy to see that a finite rated extremal system of sets from
Definition~\ref{Def:RES-fin} is a particular case of
Definition~\ref{Def:RES-Inf}. Indeed, suppose that $\ox$ is a
rated extremal point of rank $\al\in[0,1)$ for a finite set system
$\{\O_1,\ldots,\O_m\}$, i.e., condition \eqref{eq:RES-fin} is
satisfied. Defining $R(r):=\frac{\gg}{r^{1-\al}}$, we have that
$rR(r)\to 0$ and $R(r)\to\infty$ as $r\to 0$; thus $R(\cdot)$ is a
rate function while condition \eqref{eq:RES-Inf} is
satisfied.\vspace*{0.05in}

Let us discuss some specific features of the rated extremality in
Definition~\ref{Def:RES-Inf} for the case of infinite systems. For
simplicity we denote $R=R(r)$ in what follows if no confusion
arises.

\begin{Remark}[Growth condition in rated extremality]
Observe that, although $\{\O_i\}_{i\in T}$ is an infinite system
in Definition~\ref{Def:RES-Inf}, the rated extremality therein
involves only {\em finitely many} sets for each given accuracy
$\ve>0$. The imposed requirement $|I|^{3/2}=o(R)$ guarantees that
$|I|^{3/2}$ {\em grows slower} than $R$, which is very crucial in
our proof of the extremal principle below. In other words, the
number of sets involved must not be {\em too large}; otherwise the
result is trivial. We prove in Theorem~\ref{Thm:REP-Inf} that the
rate $|I|^{3/2}=o(R)$ ensures the validity of the rated extremal
principle, where the number $r$ measures {\em how far} the sets
are shifted.
\end{Remark}

Define next extremality conditions for infinite systems of sets,
which we are going to justify as an appropriate extremal principle
in what follows. These conditions are of the approximate extremal
principle type expressed in terms of of Fr\'echet normals at
points nearby the reference one.

\begin{Definition}[Rated extremality conditions for infinite systems]
\label{Def:REP-Inf} Let $\{\O_i\}_{i\in T}$ be a system of
nonempty subsets of $X$ indexed by an arbitrary set $T$, and let
$\ox\in\bigcap_{t\in T}\O_i$. We say that the set system
$\{\O_i\}_{i\in T}$ satisfies the {\sc rated extremal principle}
at $\ox$ if for any $\ve>0$ there exist a number $r\in(0,\ve)$, an
finite index subset $I\subset T$ with cardinality $|I|r<\ve$,
points $x_i\in\O_i\cap B(\ox,\ve)$, and dual elements
$x^*_i\in\Hat N(x_i;\O_i)+r\B^*$ for $i\in I$ such that
\begin{equation}\label{eq:REP-Inf}
\sum_{i\in I}x^*_i=0\;\mbox{ and }\;\sum_{i\in I}\|x^*_i\|^2=1.
\end{equation}
\end{Definition}

Observe that when a system consists of finitely many sets
$\{\O_1,\ldots,\O_m\}$ with $|I|=m$, we put the other sets equal
to the whole space $X$ and reduce Definition~\ref{Def:RES-Inf} in
this case to the conventional conditions of the approximate
extremal principle for finite systems of sets; see
Section~2.\vspace*{0.05in}

Now we address the nontriviality issue for the introduced version
of the extremal principle for infinite set systems. It is
appropriate to say (roughly speaking) that a version of the
extremal principle is {\em trivial} if all the information is
obtained from only one set of the system while the other sets
contribute nothing; i.e., if $y^*_i=0\in\Hat N(x_i;\O_i)$ for all
but one index $i$. This issue was first addressed in
\cite{MorPh10a}, where it has been shown that a ``natural"
extension of the approximate extremal principle for countable
systems is trivial.\vspace*{0.05in}

The next proposition justifies the nontriviality of the rated
extremal principle for infinite set systems proposed in
Definition~\ref{Def:REP-Inf}.

\begin{Proposition}[Nontriviality of rated extremality conditions for
infinite systems] Let $\{\O_i\}_{i\in T}$ be a system of set
satisfying the extremality conditions of
Definition~{\rm\ref{Def:REP-Inf}} at some point
$\ox\in\bigcap_{t\in T}\O_i$. Then the rated extremal principle
defined by these conditions is nontrivial.
\end{Proposition}
{\bf Proof.} Suppose on the contrary that the rated extremal
principle of Definition~\ref{Def:REP-Inf} is trivial, i.e., there
is $i_0\in T$ (say $i_0=1$) and $y^*_i\in X^*$ as $i\in T$ such
that
\begin{equation*}
x^*_i\in y^*_i+r\B^*\subset\Hat N(x_i;\O_i)+r\B^*\;\mbox{ for all
}\;i\in I,
\end{equation*}
\begin{equation*}
\sum_{i\in I}x^*_i=0,\;\;\sum_{i\in I}\|x^*_i\|^2=1,\;\mbox{ and
}\;y^*_i=0\;\mbox{ whenever }\;i\in I\setminus\{1\}
\end{equation*}
in the notation of Definition~\ref{Def:RES-Inf}. It follows that
$\|x^*_i\|\le r$ for all $i\in I\setminus\{1\}$ implying that
\begin{equation*}
\Big\|y^*_1+\sum_{i\ne1}x^*_i\Big\|\le r\;\mbox{ and
}\;\|y^*_1\|\le|I|r.
\end{equation*}
Thus we arrive at the relationships
\begin{equation*}
\sum_{i\in I}\|x^*_i\|^2<(\|y^*_1\|+r)^2+\sum_{i\ne 1}r^2\le
|I|^2r^2+2|I|r^2+r^2+(|I|-1)r^2<C\ve^2\dn 0
\end{equation*}
as $\ve\dn 0$, a contradiction. This justifies the nontriviality
of the rated extremal principle. $\h$\vspace*{0.05in}

Observe further that the extremal principle of
Definition~\ref{Def:REP-Inf} may be trivial is the rate condition
$|I|r<\ve$ is not imposed. The following example describes a
general setting when this happens.

\begin{Example}[The rate condition is essential for nontriviality]
Assume that the condition $|I|r<\ve$ is violated in the framework
of Definition~\ref{Def:REP-Inf}. Fix $\nu>0$, suppose that
$I=\{1,\ldots,N\}$ with $Nr>\nu$, pick some $u^*\in\Hat
N(x_1;\O_1)$ with the norm $\|u^*\|=\nu$, and define the dual
elements
\begin{equation*}
\begin{aligned}
&x^*_1:=u^*-\frac{u^*}{N}\in\Hat N(x_1;\O_1)+r\B^*,\\
&x^*_i:=0-\frac{u^*}{N}\in\Hat N(x_i;\O_i)+r\B^*\;\mbox{ for all
}\;i=2,\ldots,N.
\end{aligned}
\end{equation*}
Then we have the relationships
\begin{equation*}
x^*_1+\ldots+x^*_N=0\;\mbox{ and
}\;\|x^*_1\|^2+\ldots+\|x^*_N\|^2>\frac{\nu^2}{4},
\end{equation*}
which imply the triviality of the rated extremal principle by
rescaling.
\end{Example}

Now we are ready to derive the main result of this section, which
justifies the validity of the rated extremal principle for rated
extremal points of infinite systems of closed sets in Asplund
spaces.

\begin{Theorem}[Rated extremal principle for infinite systems]
\label{Thm:REP-Inf} Let $\{\O_i\}_{i\in T}$ be a system of closed
sets in an Asplund space $X$, and let $\ox$ be a rated extremal
point of this system. Then the rated extremality conditions of
Definition~{\rm\ref{Def:REP-Inf}} are satisfied for
$\{\O_i\}_{i\in T}$ at $\ox$.
\end{Theorem}
{\bf Proof.} Given $\ve>0$, take $r=\sup_i\|a_i\|$ sufficiently
small and pick the corresponding index subset $I=\{1,\ldots,N\}$
with $N^{3/2}=o(R)$ from Definition~\ref{Def:RES-Inf}. Consider
the product space $X^N$ with the norm of $z=(x_1,\ldots,x_N)\in
X^N$ given by
\begin{equation*}
\|z\|:=(\|x_1\|^2+\ldots+\|x_N\|^2)^\frac{1}{2}
\end{equation*}
and define a function $\ph\colon X^N\to\R$ by
\begin{equation}\label{ph1}
\ph(z):=\left(\sum_{i=2}^N\|(x_1-a_1)-(x_i-a_i)\|^2\right)^\frac{1}{2}.
\end{equation}
To proceed, denote $\oz:=(\ox,\ox,\ldots,\ox)\in
\O_1\times\ldots\times\O_N$ and form the set
\begin{equation}
\label{eq:REP-Inf W}
W:=\Big(\O_1\times\ldots\times\O_N\Big)\cap\Big(B\big(\ox,(R-1)r\big)\times\ldots\times
B\big(\ox,(R-1)r\big)\Big),
\end{equation}
which is nonempty and closed. We conclude that $\ph(z)>0$ for all
$z\in W$. Indeed, suppose on the contrary that $\ph(z)=0$ for some
$z=(x_1,\ldots,x_N)\in W$ and get by the estimates
$\|x_1-a_1-\ox\|\le\|x_1-\ox\|+\|a_1\|\le(R-1)r+r=Rr$ the
relationships
\begin{equation*}
x_1-a_1=\ldots=x_N-a_N\in\bigcap_{i=1}^N(\O_i-a_i)\cap
B(\ox,Rr)\ne\emp,
\end{equation*}
which contradict the extremality condition \eqref{eq:RES-Inf}.
Observe further that
\begin{equation*}
\ph(\oz)=\left(\sum_{i=2}^N\|a_1-a_i\|^2\right)^\frac{1}{2}<
2r\sqrt{N}\le\inf_{z\in W}\ph(z)+2rN^\frac{1}{2}.
\end{equation*}
Now we apply Ekeland's variational principle (see, e.g.,
\cite[Theorem~2.26]{m-book1}) with the parameters
\begin{equation*}
\ve:=2rN^\frac{1}{2}\;\mbox{ and
}\;\lm:=rR^\frac{1}{2}N^\frac{3}{4}
\end{equation*}
to the lower semicontinuous and bounded from below function
$\ph(z)+\dd(z;W)$ on $X^N$ and find in this way $z_0\in W$ such
that $\|z_0-\oz\|\le\lm$ and that $z_0$ minimizes the perturbed
function
\begin{equation}\label{beta}
\ph(z)+\beta\|z-z_0\|+\dd(z;W)\;\mbox{ on }\;z\in X^N\;\mbox{ with
}\;\be:=\frac{\ve}{\lm}=\frac{2}{R^\frac{1}{2}N^\frac{1}{4}}.
\end{equation}
By the imposed growth condition $N^\frac{3}{2}=o(R)$ as $r\dn 0$
we have
\begin{equation*}
\begin{aligned}
&\ve=2rN^\frac{1}{2}=r\cdot o(R^\frac{1}{3})\le r\cdot
o\Big(\frac{1}{r}\Big)^
\frac{1}{3}\le r\cdot o\Big(\frac{1}{r}\Big)\to 0,\\
&\frac{\lm}{Rr}=\frac{rR^\frac{1}{2}N^\frac{3}{4}}{Rr}=\frac{N^\frac{3}{4}}
{R^\frac{1}{2}}\to 0,\\
&N\be=\frac{2N}{R^\frac{1}{2}N^\frac{1}{4}}=\frac{2N^\frac{3}{4}}{R^\frac{1}{2}}=2
\Big(\frac{N^\frac{3}{2}}{R}\Big)^\frac{1}{2}\to 0\;\mbox{ as
}\;r\dn 0.
\end{aligned}
\end{equation*}
Thus $\lm=o(Rr)$ and $\be\dn 0$ as $r\dn 0$ for the quantity $\be$
defined in \eqref{beta}. Taking into account that the function
$\ph(\cdot)+\be\|\cdot-z_0\|$ is obviously Lipschitz continuous
around $\oz$, we apply to this sum the subdifferential fuzzy sum
rule from \cite[Lemma~2.32]{m-book1}. This allows us to find, for
any given number $\eta>0$, elements $z_1=(y_1,\ldots,y_N)\in
z_0+\eta\B$ and $z_2=(x_1,\ldots,x_N)\in z_0+\eta\B$ such that
\begin{equation}
\label{eq:REP-Inf1}
\big|\ph(z_1)+\be\|z_1-z_0\|-\ph(z_0)\big|\le\eta,\;z_2\in
W,\;\mbox{ and}
\end{equation}
\begin{equation}
\label{eq:REP-Inf2}
0\in\Hat\partial\Big(\ph(\cdot)+\be\|\cdot-z_0\|\Big)(z_1)+\Hat
N(z_2;W)+\eta\B^*.
\end{equation}

Our next step is to explore formula (\ref{eq:REP-Inf2}). Since
$\ph(z_0)>0$, we choose
\begin{equation*}
\eta\le\min\Big\{\be,\lm,\frac{\ph(z_0)}{2(1+\be)}\Big\}.
\end{equation*}
Then it follows from (\ref{eq:REP-Inf1}) that
\begin{equation*}
|\ph(z_1)-\ph(z_0)|\le(1+\be)\eta\le
(1+\be)\frac{\ph(z_0)}{2(1+\be)}=\frac{\ph(z_0)}{2},
\end{equation*}
which implies that $\ph(z_1)=:\al>0$. It is easy to see that the
function $\ph(\cdot)$ in \eqref{ph1} is convex. Applying the
Moreau-Rockafellar theorem of convex analysis gives us
\begin{equation}
\label{eq:REP-Inf3}
\Hat\partial\Big(\ph(\cdot)+\be\|\cdot-z_0\|\Big)(z_1)=\Hat\partial\ph(z_1)+\be\B^*,
\end{equation}
where the Fr\'echet subdifferentials on both sides of
\eqref{eq:REP-Inf3} reduce to the classical subdifferential of
convex functions. By the structure of $\ph$ in \eqref{ph1} and
that of $z_1$ we have
\begin{equation*}
\ph(z_1)=\left(\sum_{i=2}^N\|(y_1-a_1)-(y_i-a_i)\|^2\right)^\frac{1}{2}.
\end{equation*}
Denote further $\xi_i:=y_1-a_1-y_i+a_i$ for $i=2,\ldots,N$ and
observe that
$\al=\ph(z_1)=\Big(\sum_{i=2}^N\|\xi_i\|^2\Big)^\frac{1}{2}$.
Since the square root function is smooth at nonzero point, we
apply the chain rule of convex analysis to derive that any element
$(y^*_1,\ldots,y^*_N)\in\Hat\partial\ph(z_1)$ has the
representation
\begin{equation*}
y^*_i=\left\{
\begin{aligned}
&-\frac{u^*_i}{\al}\cdot\|\xi_i\|&{\rm if}\ \xi_i\ne 0,\\
&0&{\rm if}\ \xi_i=0,
\end{aligned}\qquad i=2,\ldots,N,
\right.
\end{equation*}
and $y^*_1=-y^*_2-y^*_3-\ldots-y^*_N$, where
$u^*_i\in\Hat\partial\|\cdot\|(\xi_i)$ is a subgradient of the
norm function calculated at the nonzero point $\xi_i$; hence
$\|u^*_i\|=1$. This yields that
\begin{equation*}
\|y^*_2\|^2+\ldots+\|y^*_N\|^2=1\;\mbox{ and }\;
\|y^*_1\|^2+\ldots+\|y^*_N\|^2\ge 1.
\end{equation*}
On the other hand, we have the estimates
\begin{equation*}
\|z_2-\oz\|\leq|z_2-z_0\|+\|z_0-\oz\|\le\eta+\lm\le 2\lm=o(Rr)
\end{equation*}
for $z_2=(x_1,\ldots,x_N)$ and hence
$\|x_i-\ox\|<\|z_2-\oz\|=o(Rr)$ for $i=1,\ldots,N$. The latter
ensures that each component $x_i$ lies in the interior of the ball
$B(\ox,(R-1)r)$. Furthermore, it follows from the structure of $W$
in (\ref{eq:REP-Inf W}) and the product formula for Fr\'echet
normals that
\begin{equation*}
\Hat N(z_2;W)=\Hat N\big(z_2;\O_1\times\ldots\times\O_N\big)=\Hat
N(x_1;\O_1)\times\ldots\times\Hat N(x_N;\O_N),
\end{equation*}
which implies by combining with (\ref{eq:REP-Inf2}) and
(\ref{eq:REP-Inf3}) the existence of
$(y^*_1,\ldots,y^*_N)\in\Hat\partial\ph(z_1)$ satisfying
\begin{align*}
&0\in y^*_i+\Hat N(x_i;\O_i)+2\be\B^*,\;\;\|x_i-\ox\|<2\lm\to 0\;\mbox{ as }\;r\dn 0,\\
&y^*_1+\ldots+y^*_N=0,\;\mbox{ and }\;
\|y^*_1\|^2+\ldots+\|y^*_N\|^2>1.
\end{align*}
Finally, replace $y^*_i$ by $-y^*_i$ and get from the above that
\begin{equation*}
\begin{aligned}
&y^*_i\in\Hat N(x_i;\O_i)+2\be\B^*,\;\;\|x_i-\ox\|<2\lm\to 0,\\
&\mbox{for }\;i=1,\ldots,N,\;\;N\be\to 0\;\mbox{ as }\;r\dn 0,\\
&y^*_1+\ldots+y^*_N=0,\;\mbox{ and }\;
\|y^*_1\|^2+\ldots+\|y^*_N\|^2\ge 1,
\end{aligned}
\end{equation*}
which gives all the relationships of the rated extremal principle
and completes the proof of the theorem. $\h$\vspace*{0.05in}

From the proof above we can distill some quantitative estimates
for the elements involved in the relationships of the rated
extremal principle.

\begin{Remark}[Quantitative estimates in the rated extremal
principle]\label{Rem:REP-inf} The proof of
Theorem~\ref{Thm:REP-Inf} essentially uses the growth assumptions
$N^{3/2}=o(R)$ and $R\le\tfrac{M}{r}$ on rated extremal points.
Observe in fact that the given proof allows us to make the
following {\em quantitative conclusions}: For any $\ve>0$ there
exist a number $r\in(0,\ve)$, an index subset
$I=\{j_1,\ldots,j_N\}$ with $N^{3/2}=o(R(r))$, and elements
\begin{equation*}
y^*_i\in\Hat N(x_i;\O_i)\;\mbox{ with }\;\|x_i-\ox\|\le
2rR^\frac{1}{2}N^\frac{3}{4}\;\mbox{ for all }\;i\in I
\end{equation*}
satisfying the relationships
\begin{equation*}
\|y^*_{j_1}+\ldots+y^*_{j_N}\|\le
2N\be=\frac{4N^\frac{3}{4}}{R^\frac{1}{2}}\ \mbox{ and }\
\|y^*_{j_1}\|^2+\ldots+\|y^*_{j_N}\|^2\ge 1.
\end{equation*}
Similar but somewhat different quantitative statement can be also
made: For any rated extremal point $\ox$ of the system
$\{\O_i\}_{i\in T}$ with a rate function $R(r)=O(r)$ there is a
constant $C>0$ such that whenever $\ve>0$ there exist a number
$r\in(0,\ve)$, an index subset $I=\{j_1,\ldots,j_N\}$ with
$N^{3/2}=o(\frac{1}{r})$, and elements
\begin{equation*}
y^*_i\in\Hat N(x_i;\O_i)\;\mbox{ with }\;\|x_i-\ox\|\le
C\sqrt{rN^\frac{3}{2}}\;\mbox{ for all }\;i\in I
\end{equation*}
satisfying the estimates
\begin{equation*}
\|y^*_{j_1}+\ldots+y^*_{j_N}\|\le C\sqrt{rN^\frac{3}{2}}\ \mbox{
and }\ \|y^*_{j_1}\|^2+\ldots+\|y^*_{j_N}\|^2\ge 1.
\end{equation*}
\end{Remark}\vspace*{0.1in}

In the last part of this section we introduce and study a certain
notion of {\em perturbed extremality} for arbitrary (finite or
infinite) set systems and compare it, in particular, with the
notion of linear subextremality known for systems of two sets.
Given two sets $\O_1,\O_2\subset X$, the number
\begin{equation*}
\vartheta(\O_1,\O_2):=\sup\big\{\nu\ge 0\big|\;\nu\B\subset\O_1-\O_2\big\}
\end{equation*}
is known as the {\em measure of overlapping} for these sets
\cite{Kruger06}. We say that the system $\{\O_1,\O_2\}$ is {\em
linear subextremal} \cite[Subsection~5.4.1]{m-book2} around $\ox$
if
\begin{equation}\label{eq:lextr}
\vartheta_{lin}(\O_1,\O_2,\ox):=
\liminf_{x_1\stackrel{\O_1}{\to}\ox,x_2\stackrel{\O_2}{\to}\ox\atop
r\downarrow 0} \frac{\vartheta\Big([\O_1-x_1]\cap
r\B,[\O_2-x_2]\cap r\B\Big)}{r}=0,
\end{equation}
which is called ``weak stationarity" in \cite{Kruger06}; see
\cite{Kruger06,m-book2} for more discussions and references. It is
proved in \cite{Kruger06} and \cite[Theorem~5.88]{m-book2} that
the linear subextremality of a closed set system $\{\O_1,\O_2\}$
around $\ox$ is equivalent, in the Asplund space setting, to the
validity of the approximate extremal principle for $\{\O_1,\O_2\}$
at $\ox$.

Our goal in what follows is to define a perturbed version of rated
extremality, which is applied to infinite set systems while
extends linear subextremality for systems of two sets as well.
Given an $R$-rated extremal system of sets $\{\O_i\}_{i\in T}$
from Definition~\ref{Def:RES-Inf}, we get that for any $\ve>0$
there are $r=\sup\|a_i\|$, $R=R(r)$, and $I\subset T$ satisfying
\begin{equation}\label{eq:RES-Inf2}
\bigcap_{i\in I}\big(\O_i-\ox-a_i\big)\cap (rR)\B=\emp.
\end{equation}
Let us now perturb (\ref{eq:RES-Inf2}) by replacing $\ox$ with
some $x_i\in\O_i\cap B_\ve(\ox)$ and arrive at the following
construction.

\begin{Definition}[Perturbed extremal systems]\label{Def:PES} Let
$\{\O_i\}_{i\in T}$ be a system of nonempty sets in $X$, and let
$\ox\in\bigcap_{i\in T}\O_i$. We say that $\ox$ is {\sc
$R$-perturbed extremal point} of $\{\O_i,i\in T\}$ if for any
$\ve>0$ there exist $r=\sup_{i\in I}\|a_i\|<\ve$, $I\subset T$
with $|I|^{3/2}=o(R)$, and $x_i\in\O_i\cap B_\ve(\ox)$ as $i\in I$
such that
\begin{equation}\label{pert}
\bigcap_{i\in I}\big(\O_i-x_i-a_i\big)\cap (rR)\B=\emp.
\end{equation}
In this case we say that $\{\O_i\}_{i\in T}$ is an {\sc
$R$-perturbed extremal system} at $\ox$.
\end{Definition}

The next proposition establishes a connection between linear
subextremality and perturbed extremality for systems of two sets
$\{\O_1,\O_2\}$.

\begin{Proposition}[Perturbed extremality from linear subextremality] Let a set system
$\{\O_1,\O_2,\ox\}$ be linearly subextremal around $\ox$. Then it
is an $R$-perturbed extremal system at this point.
\end{Proposition}

{\bf Proof.} Employing the definition of linear subextremality,
for any $\ve>0$ sufficiently small we find $x_i\in\O_i\cap
B_\ve(\ox)$ and $r'<\ve$ such that
\begin{equation*}
\vt\big([\O_1-x_1]\cap r'\B,[\O_2-x_2]\cap r'\B\big)<r'\ve.
\end{equation*}
This implies the existence of a vector $a\in X$ satisfying
$\|a\|\le r'\ve$ and
\begin{equation*}
a\not\in \Big([\O_1-x_1]\cap r'\B\Big)-\Big([\O_2-x_2]\cap
r'\B\Big),
\end{equation*}
which ensures in turn that
\begin{equation}\label{s1}
\Big([\O_1-x_1]\cap r'\B-\frac{a}{2}\Big)\cap\Big([\O_2-x_2]\cap
r'\B+\frac{a}{2}\Big)=\emp.
\end{equation}
Let us show that the latter implies the fulfillment of
\begin{equation}\label{s2}
\Big[\O_1-x_1-\frac{a}{2}\Big]\cap\Big[\O_2-x_2+\frac{a}{2}\Big]\cap
\frac{r'}{2}\B=\emp.
\end{equation}
Indeed, suppose that \eqref{s2} does not hold and pick $\xi\in X$
from the left-hand side set in \eqref{s2}. Since
$\xi+\tfrac{a}{2}\in\O_1-x_1$ and $\|\xi\|\le\tfrac{r'}{2}$, we
have
\begin{equation*}
\Big\|\xi_\frac{a}{2}\Big\|\le\frac{r'}{2}+\frac{r'\ve}{2}\le\frac{r'}{2}
+\frac{r'}{2}=r'
\end{equation*}
and consequently $\xi\in[\O_1-x_1]\cap r'\B-\disp\frac{a}{2}$.
Similarly we get $\xi\in[\O_2-x_2]\cap r'\B-\disp\frac{a}{2}$.
This clearly contradicts \eqref{s1} and thus justifies the claimed
relationship \eqref{s2}.

By setting $r:=\disp\frac{\|a\|}{2}$, out remaining task is to
construct a continuous function $\colon\R_+\to\R_+$ such that
$R(r)\to\infty$ as $r\dn 0$ and that for each $\ve>0$ there is
$r<\ve$ satisfying
\begin{equation*}
\disp\Big[\O_1-x_1-\frac{a}{2}\Big]\cap\Big[\O_2-x_2+\frac{a}{2}\Big]\cap
(rR)\B=\emp.
\end{equation*}
We first construct such a function along a sequence $r_k\dn 0$ as
$k\to\infty$. Picking $\ve_k\dn 0$, find $r'_k<\ve_k$ and select
$a_k\in X$ with $\|a_k\|\le r'_k\ve_k$ such that the sequence of
$\|a_k\|$ is decreasing. Then define $r_k:=\disp\frac{\|a_k\|}{2}$
and $R(\ve_k):=\disp\frac{1}{\ve_k}$. It follows from the
constructions above that
\begin{equation*}
r_kR(r_k)\le r'_k\ve_k\frac{1}{\ve_k}=r'_k,\quad k\in\N.
\end{equation*}
We clearly see that the sequence $\{R(r_k)\}$ is increasing as
$r_k\dn 0$. Extending $R(\cdot)$ piecewise linearly to $\R_+$
brings us to the framework of Definition~\ref{Def:PES} and thus
completes the proof of the proposition. $\h$ \vspace{0.05in}

Finally in this section, we show the rated extremality conditions
of Definition~\ref{Def:REP-Inf} holds for $R$-perturbed extremal
points of infinite set systems from Definition~\ref{Def:PES}.

\begin{Theorem}[Rated Extremal Principle for Perturbed Systems]
\label{Thm:REP-Per} Let $\ox$ be an $R$-perturbed extremal point
of a closed set system $\{\O_i\}_{i\in T}$ in an Asplund space
$X$. Then the rated extremal principle holds for this system at
$\ox$.
\end{Theorem}

{\bf Proof.} Fix $\ve>0$ and find $I$, $\{x_i\}_{i\in I}$, and
$\{a_i\}_{i\in I}$ from Definition~\ref{Def:PES} such that
\begin{equation*}
\bigcap_{i\in I}\big(\O_i-x_i-a_i\big)\cap(rR)\B=\emp.
\end{equation*}
For convenience denote $I:=\{1,\ldots,N\}$ and define
$$
\O:=\Big\{(u_1,\ldots,u_N)\in X^N\,\Big|\,u_i\in\O_i\cap(x_i+rR\B),\;i\in
I\Big\}.
$$
For any $z=(u_1,\ldots,u_N)\in\O$ consider the function
\begin{equation*}
\ph(z):=\Big(\sum_{i=2}^N\|(u_1-x_1-a_1)-(u_i-x_i-a_i)\|^2\Big)^
\frac{1}{2}>0.
\end{equation*}
Furthermore, for $\oz=(x_1,\ldots,x_N)$ we have the estimates
\begin{equation*}
\disp\ph(\oz)=\Big(\sum_{i=2}^N\|a_1-a_i\|^2\Big)^\frac{1}{2}<
2r\sqrt{N}\le\inf_{z\in \O}\ph(z)+2rN^\frac{1}{2}.
\end{equation*}
The rest of the proof follows the arguments in the proof of
Theorem~\ref{Thm:REP-Inf}.$\h$

\section{Calculus Rules for Rated Normals to Infinite Intersections}
\setcounter{equation}{0}

In the concluding section of the paper we apply the rated extremal
principle of Section~4 to deriving some calculus rules for general
normals to infinite set intersections, which are closely related
to necessary optimality conditions in problems of semi-infinite
and infinite programming. Unless otherwise stated, the spaces
below are Asplund and the sets under consideration are closed
around reference points. As in Section~4, we often drop the
subscript ``$r$'' for simplicity in the notation of rate functions
$R_r=R(r)$ if no confusion arises. In addition, we always assume
that rate functions are continuous.\vspace*{0.05in}

We start with the following definition of {\em rated normals} to
set intersections.

\begin{Definition}[Rated normals to set intersection]\label{Def:Rnor}
Let $\O:=\bigcap_{i\in T}\O_i$, and let $\ox\in\O$. We say that a
dual element $x^*\in X^*$ is an {\sc $R$-normal} to the set
intersection $\O$ if for any $r\dn 0$ there is $I=I(r)\subset T$
of cardinality $|I|^{3/2}=o(R_r)$ such that
\begin{equation}\label{eq:Rnor}
\la x^*,x-\ox\ra-r\|x-\ox\|<r\ \mbox{ for all }\ x\in\bigcap_{i\in
I}\O_i\cap B(\ox,rR_r).
\end{equation}
\end{Definition}
The next proposition reveals relationships between \F\ and
$R$-normals to set intersections.

\begin{Proposition}[Rated normals versus \F\ normals to set intersections]
\label{Prop:R-Fnor} Let $\ox\in\O=\bigcap_{i\in I}\O_i$. Then any
$R$-normal to $\O$ at $\ox$ is a \F\ normal to $\O$ at $\ox$. The
converse holds if $I$ is finite.
\end{Proposition}
{\bf Proof.} Assume $x^*$ is an $R$-normal to $\O$ at $\ox$ with
some rate function $R(r)$ while $x^*$ is not a \F\ normal to $\O$
at this point. Hence there are $\dd>0$ and a sequence
$x_k\st{\O}{\to}\ox$ such that $\dd\|x_k-\ox\|<\la x^*,x_k-\ox\ra$
for all $k\in\N$. Hence $x_k\ne\ox$ and
\begin{equation*}
\dd\|x_k-\ox\|<\la x^*,x_k-\ox\ra<r\|x_k-\ox\|+r
\end{equation*}
whenever $\|x_k-\ox\|\le rR$. Now suppose that $rR=M>0$ for some
$M$ and then fix a number $k\in\N$ such that $\|x_k-\ox\|\le rR$.
Letting $r\dn0$, we arrive at the contradiction
$\dd\|x_k-\ox\|\leq 0$.

Consider next the remaining case when $rR\to 0$ as $r\dn 0$ and find $r_k>0$
sufficiently small so that $\|x_k-\ox\|=r_kR(r_k)$ due to the continuity of $R$ and the
convergence $rR\st{r\dn0}{\lto}0$. It follows that
\begin{equation*}
\dd r_kR(r_k)< r^2_kR(r_k)+r_k\;\mbox{ and hence }\;\dd<
r_k+\frac{1}{R(r_k)}, \quad k\in\N,
\end{equation*}
which gives a contradiction as $k\to\infty$. Thus $x^*$ is a \F\
normal to $\O$ at $\ox$.

Conversely, assume that the index set $I$ is finite, i.e.,
$I=\{1,\ldots,N\}$, and that $x^*$ is a \F\ normal. Then for any
$r>0$ we have by \eqref{eq:e-nor} that
\begin{equation*}
\la x^*,x-\ox\ra-r\|x-\ox\|\le 0\ \mbox{ for all }\
x\in\bigcap_{i=1}^N\O_i\cap U,
\end{equation*}
where $U$ is a neighborhood of $\ox$. This clearly implies
(\ref{eq:Rnor}) with any rate function $R$, which ensures that
$x^*$ is an $R$-normal to $\O$ at $\ox$ and thus completes the
proof of the proposition. $\h$\vspace*{0.05in}

The next example concerns infinite systems of convex sets in
$\R^2$. It illustrates the way of computing $R$-normals to
infinite intersections and shows that $R$-normals in this case
reduce to usual ones.

\begin{Example}[Rated normals for infinite systems]
\label{Ex:Rnor-inf} Let $m\ge 4$ be a fixed integer. Consider an
infinite system of convex sets $\{\O_k\}_{k\in\N}$ in $\R^2$
defined as the epigraphs of the convex and smooth functions
$$
g_k(x):=
\begin{cases}k^mx^2 &\mbox{ for }\;x\ge 0,\\ 0 &\mbox{ for }\;x<0,
\end{cases}\quad k=1,2,\ldots.
$$
Let $\ox:=(0,0)$, $\O:=\bigcap_{k=1}^\infty\O_k$, and let
$R=R(r)=r^{\al-1}$ for some $\al\in(0,\tfrac{2}{11})$. We
obviously get $\O=\R_-\times\R_+$ and $N(\ox;\O)=\R_+\times\R_-$.
Let us verify that $x^*=(1,0)$ is an $R$-normal to $\O$ at $\ox$,
which implies the whole normal cone $N(\ox;\O)$ consists of
$R$-normals.

To proceed, fix any $r>0$ sufficiently small and denote by $k_0$
the {\em smallest} integer such that
\begin{equation*}
\max\Big\{\frac{1}{4
r^2},\frac{1}{4r^{2+\al}}\Big\}=\frac{1}{4r^{2+\al}}\le k^m_0.
\end{equation*}
Now consider $I:=\{1,\ldots,k_0\}$ and check that
\begin{equation*}
k_0\le
\Big(\frac{1}{4r^{2+\al}}\Big)^{1/m}+1<\frac{1}{r^\frac{2+\al}{m}}.
\end{equation*}
Since $1-\frac{3}{2m}(2+\al)-\al\ge 1-\frac{3}{8}(2+\al)-\al\ge
\frac{1}{4}-\frac{11}{8}\al>0$, it follows that
\begin{equation*}
\frac{|I|^{3/2}}{R}<\frac{r^{1-\al}}{r^{\frac{3(2+\al)}{2m}}}=r^{1-
\frac{3}{2m}(2+\al)-\al}\to 0\quad\mbox{when }\ r\dn 0.
\end{equation*}
Defining further $\O_0:=\bigcap_{k=1}^{k_0}\O_k$, it remains to
show that
\begin{equation}\label{k0}
\la x^*,x\ra -r\|x\|<r\ \mbox{ for all }\ x\in\O_0\cap B(0;r
R).
\end{equation}
To verify \eqref{k0}, take $x:=(t,s)$ and consider only the case
when $t>0$, since the other case of $t\le 0$ is obvious. For $t>0$
we have $s\ge k^m_0 t^2$ and
\begin{equation}\label{k1}
\la x^*,x\ra -r\|x\|=t-r\sqrt{t^2+s^2}\le t\Big(1
-r\sqrt{1+k^{2m}_0 t^2}\Big)<t\big(1-r k^m_0 t\big)=-r k^m_0
t^2+t=:f(t).
\end{equation}
It follows from $\|x\|\le r R=r^\al$ that
\begin{equation*}
r^\al\ge\sqrt{t^2+s^2}\ge t\sqrt{1+k^{2m}_0t^2}>k^m_0 t^2
\end{equation*}
and hence $t<\big(\frac{r^{\al}}{K^m}\big)^{1/2}$. The latter
implies that for all $x=(t,s)\in\O_0\cap B(0;r R)$ with $t>0$ we
have
\begin{equation*}
\la x^*,x\ra -r\|x\|<f(t)\le\sup_{[0,a]}f(t)\quad\mbox{with }\
a:=\Big(\frac{r^{\al}}{k^m_0}\Big)^{1/2}\ge\frac{1}{2r k^m_0}.
\end{equation*}
Observe finally that the function $f(t)$ in \eqref{k1} attains its
maximum on [0,a] at the point $t=\tfrac{1}{2r k^m_0}$ and that
\begin{equation*}
\disp\sup_{[0,a]}f(t)=-r k_0\frac{1}{4r^2 k^{2m}_0}+\frac{1}{2r
k^m_0}=\disp\frac{1}{4 r k^m_0}\le r.
\end{equation*}
Combining all the above, we arrive at \eqref{k0} and thus achieve
our goals in this example.
\end{Example}\vspace*{0.1in}

The next example related to the previous one involves the notion
of equicontinuity for systems of mappings. Given $f_i\colon X\to
Y$, $i\in T$, we say that the system $\{f_i\}_{i\in T}$ is {\em
equicontinuous}  at $\ox$ if for any $\ve>0$ there is $\dd>0$ such
that $\|f_i(x)-f_i(\ox)\|<\ve$ for all $x\in B(\ox,\dd)$ and $i\in
T$. This notion has been recently exploited in \cite{Seidman10} in
the framework of variational analysis; see Remark~\ref{equi}.

\begin{Example}[Non-equicontinuity of gradient and normal systems]
\label{Ex:nequi} Given an integer $m\ge 4$, define an infinite
systems of functions $\ph_k\colon\R^2\to\R$ for $k\in\N$ by
\begin{equation}\label{equi1}
\disp\ph_k(x_1,x_2):=\begin{cases}k^m x^2_1-x_2&\mbox{for }\;x_1>0,\\
-x_2&\mbox{for }\;x_1\le 0.\end{cases}
\end{equation}
It is easy to check that the system of gradients
$\{\nabla\ph_k\}_{k\in\N}$ is not equicontinuous at $\ox=(0,0)$.

Furthermore, observe that the sets $\O_k$ in
Example~\ref{Ex:Rnor-inf} can be defined by
\begin{equation}\label{ok}
\O_k:=\big\{x\in\R^2\big|\;\ph_k(x)\le 0\big\},\quad k\in\N.
\end{equation}
Given any boundary point $(x_1,x_2)$ of the set $\O_k$, we compute
the unit normal vector to $\O_k$ at $(x_1,x_2)$ by
\begin{equation*}
\xi_k(x_1,x_2)=\begin{cases}\disp\frac{1}{\sqrt{4k^{2m}x^2_1+1}}(2k^m
x_1,-1)
&\mbox{for }\;x_1>0,\\
(0,-1)&\mbox{for }\;x_1\le 0.\end{cases}
\end{equation*}
and then check the relationships for $x_1>0$:
\begin{equation*}
\|\xi_k(x_1,x_2)-\xi_k(0,0)\|^2=\frac{8k^{2m}x^2_1-2\sqrt{4k^{2m}x^2_1+1}}{4k^{2m}
x^2_1+1}\to 2\ \mbox{ as }\ k\to\infty.
\end{equation*}
The latter means that the system of $\{\xi_k\}_{k\in\N}$ is not
equicontinuous at $\ox=(0,0)$.
\end{Example}\vspace*{0.05in}

The next major result of this paper establishes a certain ``fuzzy"
intersection rule for rated normals to infinite set intersections.
Its proof is based on the rated extremal principle for infinite
set systems obtained above in Theorem~\ref{Thm:REP-Inf}. Parts of
this proof are similar to deriving a fuzzy sum rule for Fr\'echet
normals to intersections of two sets in Asplund spaces given in
\cite{Mor-Wang2002} and in \cite[Lemma~3.1]{m-book1} on the base
of the approximate extremal principle for such set systems.

\begin{Theorem}[Fuzzy intersection rule for $R$-normals]
\label{Thm:Fuz} Let $\ox\in\O:=\bigcap_{i\in T}\O_i$, and let
$x^*\in X^*$ be an $R$-normal to $\O$ at $\ox$. Then for any
$\ve>0$ there exist an index subset $I$, Fr\'echet normals
$x^*_i\in\Hat N(x_i;\O_i)$ with $\|x_i-\ox\|<\ve$ for $i\in I$,
and a number $\lm\ge 0$ such that
\begin{equation}\label{fuzzy}
\lm x^*\in\sum_{i\in I}x^*_i+\ve\B^*\;\mbox{ and
}\;\lm^2+\lm^2\|x^*\|^2+\sum_{i\in I}\|x^*_i\|^2=1.
\end{equation}
\end{Theorem}
{\bf Proof.} Without loss of generality, assume that $\ox=0$. Pick
any $x^*\in\Hat N(0;\O)$ and by Definition~\ref{Def:Rnor} for any
$r>0$ sufficiently small find an index subset $|I|^{3/2}=o(R)$
such that
\begin{equation}\label{eq:Fuz1}
\la x^*,x\ra-r\|x\|<r\ \mbox{ whenever }\ x\in\bigcap_{i\in
I}\O_i\cap (r R)\B.
\end{equation}
Then we form the following closed subsets of the Asplund space
$X\times\R$:
\begin{equation}\label{O}
\begin{aligned}
&O_1:=\Big\{(x,\al)\in X\times\R\Big|\;x\in\O_1,\;\al\le\la x^*,x\ra-r\|x\|\Big\},\\
&O_i:=\O_i\times\R_+\;\mbox{ for }\;i\in I\setminus\{1\},
\end{aligned}
\end{equation}
where $I=\{1,\ldots,N\}$  with ``$1$'' denoting the first element
of $I$ for simplicity. This leads us to
\begin{equation}\label{emp1}
\Big(O_1-(0,r)\Big)\cap\bigcap_{i\in I\setminus\{1\}}O_i\cap(r R_r)\B=\emp.
\end{equation}
Indeed, if on the contrary \eqref{emp1} does not hold, we get
$(x,\al)$ from the above intersection satisfying $\al\ge 0$,
$x\in\bigcap_{i\in I}\O_i\cap(\ve R_\ve)\B$, and
\begin{equation*}
r\leq\al+r\le\la x^*,x\ra-r\|x\|,
\end{equation*}
where the latter is due to $(x,\al+r)\in O_1$. This clearly
contradicts (\ref{eq:Fuz1}) and so justifies \eqref{emp1}. Thus we
have that $(0,0)\in X\times\R$ is a rated extremal point of the
set system $\{O_1,O_2\}$  from \eqref{O} in the sense of
Definition~\ref{Def:RES-Inf}. Applying to this system the rated
extremal principle from Theorem~\ref{Thm:REP-Inf} with taking into
account Remark~\ref{Rem:REP-inf} to find elements $(w_i,\al_i)$
and $(x^*_i\lm_i)$ for $i=1,\ldots,N$ satisfying the relationships
\begin{eqnarray}\label{O1}
\left\{\begin{array}{ll} &(x^*_i,\lm_i)\in\Hat
N\big((w_i,\al_i);O_i\big),\;\|(w_i,\al_i)\|\le 2rR^\frac{1}{2}
N^\frac{3}{4},\quad i\in I,\\
&\Big\|(x^*_1,\lm_1)+\ldots+(x^*_N,\lm_N)\Big\|\le\disp\frac{4N^\frac{3}{4}}{R^\frac{1}{2}}=
:\eta\dn 0\;\mbox{ as }\;r\dn 0,\\
&\|(x^*_1,\lm_1)\|^2+\ldots+\|(x^*_N,\lm_N)\|^2=1.
\end{array}\right.
\end{eqnarray}
By the structure of $O_i$ as $i=1,\ldots,N$ we have from the first
line of \eqref{O1} that $x^*_i\in\Hat N(w_i;\O_i)$, that $\lm_i\le
0$ for $i=2,\ldots,N$, and that
\begin{equation}
\label{eq:Fuz2} \limsup_{(x,\al)\st{O_1}{\to}(w_1,\al_1)}\frac{\la
x^*_1,x-w_1\ra+\lm_1(\al-\al_1)}{\|x-w_1\|+|\al-\al_1|}\le 0
\end{equation}
by the definition of \F\ normals. It also follows from the
structure of $O_1$ that $\lm_1\ge 0$ and
\begin{equation}
\label{eq:Fuz3}\al_1\le\la x^*,w_1\ra-r\|w_1\|.
\end{equation}
This allows us to split the situation into the follows two
cases.\\[1ex]
{\bf Case~1:} $\lm_1=0$. If inequality (\ref{eq:Fuz3}) is strict
in this case, there is a neighborhood $W$ of $w_1$ such that
\begin{equation*}
\al_1\le\la x^*,x\ra-r\|x\|\ \mbox{ for all }\ x\in \O_1\cap W.
\end{equation*}
This implies that $(x,\al_1)\in O_1$ whenever $x\in \O_1\cap W$.
Substituting $(x,\al_1)$ into (\ref{eq:Fuz2}) gives us
\begin{equation*}
\limsup_{x\st{\O_1}{\to}w_1}\frac{\la
x^*_1,x-w_1\ra}{\|x-w_1\|}\le 0,\;\mbox{ i.e., }\;x^*_1\in\Hat
N(w_1;\O_1).
\end{equation*}
If (\ref{eq:Fuz3}) holds as equality, we denote $\al:=\la x^*,x\ra
-r\|x\|$ and get
\begin{equation*}
|\al-\al_1|=\Big|\la x^*,x-w_1\ra+r(\|w_1\|-\|x\|)\Big|\le
\Big(\|x^*\|+r\Big)\|x-w_1\|,
\end{equation*}
which implies by (\ref{eq:Fuz2}) that
\begin{equation*}
\limsup_{(x,\al)\st{O_1}{\to}(w_1,\al_1)} \frac{\la
x^*_1,x-w_1\ra}{\|x-w_1\|+|\al-\al_1|}\le 0.
\end{equation*}
Thus it follows for any $\ve'>0$ sufficiently small and the number
$\al$ chosen above that
\begin{equation*}
\la x^*_1,x-w_1\ra\le\ve'\Big(\|x-w_1\|+|\al-\al_1|\Big)\le
\ve'\Big(1+\|x^*\|+r\Big)\|x-w_1\|
\end{equation*}
for all $x\in\O_1$ sufficiently closed to $w_1$. This ensures that
\begin{equation*}
\limsup_{x\st{\Lm_1}{\to}w_1}\frac{\la
x^*_1,x-w_1\ra}{\|x-w_1\|}\le 0,\;\mbox{ i.e., }\;x^*_1\in\Hat
N(w_1;\O_1)
\end{equation*}
when \eqref{eq:Fuz3} holds as equality as well as the strict
inequality. Since $\lm_1=0$ in Case~1 under consideration and
since $\lm_i\le 0$ for all $i\ge 2$, it follows that
\begin{equation*}
\lm^2_2+\ldots+\lm^2_N\le(\lm_2+\ldots+\lm_N)^2\le\eta^2.
\end{equation*}
This leads us to the estimates
\begin{equation*}
\|x^*_1\|^2+\ldots+\|x^*_N\|^2\ge 1-(\lm^2_2+\ldots+\lm^2_N)\ge
\frac{1}{2},
\end{equation*}
and thus we get from \eqref{O1} all the conclusion of the theorem
with $\lm=0$ in \eqref{fuzzy} in this case.\\[1ex]
{\bf Case~2:} $\lm_1>0$. If inequality (\ref{eq:Fuz3}) is strict
in this case, put $x:=w_1$ and get from (\ref{eq:Fuz2}) that
\begin{equation*}
\limsup_{\al\to\al_1}\frac{\lm_1(\al-\al_1)}{|\al-\al_1|}\le 0,
\end{equation*}
which yields $\lm_1=0$, a contradiction. It remains therefore to
consider the case when (\ref{eq:Fuz3}) holds as equality. Take
then a pair $(x,\al)\in O_1$ with
\begin{equation*}
x\in\O_1\setminus\{w_1\}\;\mbox{ and }\;\al=\la x^*,x\ra-r\|x\|
\end{equation*}
and hence get from (\ref{eq:Fuz3}) that
\begin{equation*}
\al-\al_1=\la x^*,x-w_1\ra+r(\|w_1\|-\|x\|),
\end{equation*}
which implies the relationships
\begin{equation*}
\la x^*_1,x-w_1\ra +\lm_1(\al-\al_1) =\la x^*_1+\lm_1 x^*,x-w_1\ra
+\lm_1 r(\|w_1\|-\|x\|),
\end{equation*}
\begin{equation*}
|\al-\al_1|\le(\|x^*\|+r)\|x-w_1\|.
\end{equation*}
On the other hand, it follows from (\ref{eq:Fuz2}) that for any
$\ve'>0$ sufficiently small there exists a neighborhood $V$ of
$w_1$ such that
\begin{equation*}
\la
x^*_1,x-w_1\ra+\lm_1(\al-\al_1)\le\lm_1\ve'r\Big(\|x-w_1\|+|\al-\al_1|\Big),
\end{equation*}
whenever $x\in\O_1\cap V$ and that
\begin{equation*}
\begin{aligned}
\la x^*_1+\lm_1 x^*,x-w_1\ra+\lm_1 r(\|w_1\|-\|x\|)&\le\lm_1
\ve'r(\|x-w_1\|+|\al-\al_1|)\\
&\le\lm_1\ve'r\Big[\|x-w_1\|+(\|x^*\|+r)\|x-w_1\|\Big]\\
&=\lm_1\ve'r\big(1+\|x^*\|+r\big)\|x-w_1\|.
\end{aligned}
\end{equation*}
Let us now choose $\ve'>0$ sufficiently small so that
\begin{equation*}
\la x^*_1+\lm_1 x^*,x-w_1\ra +\lm_1 r(\|w_1\|-\|x\|)\le\lm_1
r\|x-w_1\|.
\end{equation*}
and for all $x\in\O_1\cap V$ get the estimate
\begin{equation*}
\la x^*_1+\lm_1x^*,x-w_1\ra\le\lm_1 r\|x-w_1\|+\lm_1
r(\|x\|-\|w_1\|)\le 2\lm_1 r\|x-w_1\|.
\end{equation*}
It follows definition \eqref{eq:e-nor} of $\ve$-normals that
\begin{equation*}
x^*_1+\lm_1 x^*\in\Hat N_{2\lm_1r}(w_1;\O_1),
\end{equation*}
where $\lm_1\le 1$ by the third line of \eqref{O1}. Using the
representation of $\ve$-normals in Asplund spaces from
\cite[(2.51)]{m-book1}, we find $v\in\O_1\cap(w_1+2\lm_1 r)\B)$
such that
\begin{equation*}
x^*_1+\lm_1x^*\in \Hat N(v;\O_1)+2\lm_1 r\B^*.
\end{equation*}
Hence $\|v\|\le\|v-w_1\|+\|w_1\|\le
2\lm_1r+2rR^\frac{1}{2}N^\frac{3}{4}\le
3rR^\frac{1}{2}N^\frac{3}{4}$ and there is $\Tilde x^*_1\in\Hat
N(v;\O_1)$ with
\begin{equation*}
\lm_1x^*\in\Tilde x^*_1-x^*_1+2\lm_1 r\B^*.
\end{equation*}
Taking into account that $x^*_1+\ldots+x^*_N\in\eta\B^*$, we get
\begin{equation*}
\lm_1x^*\in\Tilde x^*_1+x^*_2+\ldots+x^*_N+(2\lm_1 r+\eta)\B^*.
\end{equation*}
On the other hand, it follows from $-x^*_1=\lm_1x^*-\Tilde
x^*_1-u^*$ with some $\|u^*\|\le 2\lm_1 r\le 2r$ that
\begin{equation*}
\|x^*_1\|^2\le\big(\lm_1\|x^*\|+\|\Tilde x^*_1\|+2r\big)^2\le
2\lm^2_1\|x^*\|^2+2\|\Tilde x^*_1\|^2+\frac{1}{4}.
\end{equation*}
Moreover, since $|\lm_1+\lm_2+\ldots+\lm_N|\le\eta\dn 0$ as $r\dn
0$ by the second line of \eqref{O1} and since $\lm_1\ge 0$ while
$\lm_i\le 0$ for $i=2,\ldots,N$, we have
\begin{equation*}
\eta^2>\lm^2_1+(\lm_2+\ldots+\lm_N)^2+2\lm_1(\lm_2+\ldots+\lm_N)>\lm^2_1+(\lm_2+\ldots
+\lm_N)^2+2\lm_1(-\lm_1-\eta)
\end{equation*}
It also follows from \eqref{O1} and $0<\lm_1<1$ that
\begin{equation*}
\lm_1^2\ge(\lm_2+\ldots\lm_N)^2-\eta^2-2\eta\lm_1\ge
\lm^2_2+\ldots+\lm^2_N-\frac{1}{4},
\end{equation*}
which leads us to the subsequent estimates
\begin{equation*}
\lm_1^2+\ldots+\lm_N^2\le 2\lm_1^2+\frac{1}{4}\;\mbox{ and}
\end{equation*}
\begin{align*}
1&\le\Big(\lm_1^2+\ldots+\lm_N^2\Big)+\Big(\|x^*_1\|^2+\ldots+\|x^*_N\|^2\Big)\\
&\le
2\lm^2_1+2\lm^2_1\|x^*\|^2+2\|\tx^*_1\|^2+\Big(\|x^*_2\|^2+\ldots+\|x^*_N\|^2\Big)
+\frac{1}{2}.
\end{align*}
This finally ensures that
\begin{equation*}
\frac{1}{4}\le\lm_1^2+\lm_1^2\|x^*\|^2+\|\Tilde
x^*_1\|^2+\|x^*_2\|^2+\ldots+\|x^*_N\|^2
\end{equation*}
and brings us to all the conclusions of the theorem with
$\lm:=\lm_1$ in \eqref{fuzzy}. $\h$

\begin{Remark}[Quantitative estimates in the intersection rule]\label{Rem:Fuz}
It can be observed directly from the proof of
Theorem~\ref{Thm:Fuz} that we get in fact the following
quantitative estimates in intersection rule obtained for infinite
set systems when $r>-$ is sufficiently small: $|I|^{3/2}=o(R)$,
\begin{equation*}
\|x_i-\ox\|<3rR^\frac{1}{2}|I|^\frac{3}{4},\;\mbox{ and }\;\lm
x^*\in\sum_{i\in
I}x^*_i+\disp\Big(2r+4\frac{|I|^\frac{3}{4}}{R^\frac{1}{2}}\Big)\B^*.
\end{equation*}
In particular, for $R=O\big(\tfrac{1}{r}\big)$, there is $C>0$
such that  all the conclusions hold with
$|I|^{3/2}=N^{3/2}=o\big(\frac{1}{r}\big)$,
\begin{equation*}
\|x_i-\ox\|<C\sqrt{rN^\frac{3}{2}},\;\mbox{ and }\;\lm
x^*\in\sum_{i\in I}x^*_i+C\sqrt{rN^\frac{3}{2}}\B^*.
\end{equation*}
\end{Remark}

\begin{Remark}[Perturbed rated normals to infinite intersections]
Inspired by our consideration of perturbed extremal systems in
Section~4, we define a perturbed version of $R$-normals to
infinite set intersections as follows: $x^*\in X^*$ is a {\em
perturbed $R$-normal} to the intersection $\O:=\bigcap_{i\in
T}\O_i$ at $\ox\in\O$ if for any $\ve>0$ there exist a number
$r>0$, an index subset $I$ with cardinality $|I|^{3/2}=o(R_r)$,
and points $x_i\in\O_i\cap B(\ox,\ve)$ as $i\in I$ such that
$r|I|<\ve$ and
\begin{equation*}
\la x^*,x\ra-r\|x\|<r\ \mbox{ whenever }\ x\in\bigcap_{i\in
I}\big(\O_i-x_i\big)\cap(r R_r)\B.
\end{equation*}
Then the corresponding version of the intersection rule from
Theorem~\ref{Thm:Fuz} can be derived for perturbed rated normals
to infinite intersections by a similar way with replacing in the
proof the rated extremal principle from Theorem~\ref{Thm:REP-Inf}
by its perturbed version from Theorem~\ref{Thm:REP-Per}.
\end{Remark}

We proceed with deriving calculus rules for the so-called {\em
limiting $R$-normals} (defined below) to infinite intersections of
sets. First we propose a new qualification conditions for infinite
systems.

\begin{Definition}[Approximate qualification condition]\label{Def:AQC}
We say that a system of sets $\{\O_i\}_{i\in T}\subset X$
satisfies the {\sc approximate qualification condition (AQC)} at
$\ox\in\bigcap_{i\in T}\O_i$ if for any $\ve\dn 0$, any finite
index subset $I_\ve\subset T$, and any Fr\'echet normals
$x^*_{i\ve}\in\Hat N(x_{i\ve};\O_i)\cap\B^*$ with
$\|x_{i\ve}-\ox\|\le\ve$ as $i\in I_\ve$ the following implication
holds:
\begin{equation}\label{eq:AQC}\Big\|\sum_{i\in I_\ve}x^*_{i\ve}\Big\|\st{\ve\dn
0}{\lto}0\Longrightarrow \sum_{i\in
I_\ve}\|x^*_{i\ve}\|^2\st{\ve\dn 0}{\lto}0.
\end{equation}
\end{Definition}

The next proposition presents verifiable conditions ensuring the
validity of AQC for finite systems of sets under the SNC property
\eqref{snc} discussed at the end of Section~3; see \cite{m-book1}
for more details.

\begin{Proposition}[AQC for finite set systems under SNC assumptions]
Let $\{\O_1,\ldots,\O_m\}$ be a finite set system satisfying the
limiting qualification condition at $\ox\in\bigcap_{i=1}^m\O_i$:
for any sequences $x_{ik}\st{\O_i}{\to}\ox$ and
$x^*_{ik}\st{w^*}{\to}x^*_i$ with $x^*_{ik}\in\Hat N(x_{ik};\O_i)$
as $k\to\infty$ and $i=1,\ldots,m$ we have
\begin{equation*}
\|x^*_{1k}+\ldots+x^*_{mk}\|\to 0\Longrightarrow
x^*_1=\ldots=x^*_m=0,
\end{equation*}
which is automatic under the normal qualification condition via
the basic normal cone \eqref{eq:Lcone}:
\begin{equation*}
\big[x^*_1+\ldots+x^*_m=0\;\mbox{ and }\;x^*_i\in
N(\ox;\O_i),\;i=1,\ldots,m\big]\Longrightarrow x^*_i=0\;\mbox{ for
all }\;i=1,\ldots,m.
\end{equation*}
Assume in addition that all but one of $\O_i$ are SNC at $\ox$.
Then the AQC is satisfied for $\{\O_1,\ldots,\O_m\}$ at $\ox$.
\end{Proposition}
{\bf Proof.} Pick $\ve_k\dn 0$, $x^*_{ik}\in\Hat
N(x_{ik};\O_i)\cap\B^*$, $\|x_{ik}-\ox\|\le\ve_k$ as
$i=1,\ldots,m$ and assume that
\begin{equation}\label{snc1}
\|x^*_{1k}+\ldots+x^*_{mk}\|\to 0\;\mbox{ as }\;k\to\infty.
\end{equation}
Taking into account that the sequences $\{x^*_{ik}\}\subset X^*$
are bounded when $X$ is Asplund, we extract from them weak$^*$
convergent subsequences and suppose with no relabeling that
$x^*_{ik}\st{w^*}{\to}x^*_i$ as $k\to\infty$ for all
$i=1,\ldots,m$. It follows from the imposed limiting qualification
condition for $\{\O_1,\ldots,\O_m\}$ at $\ox$ that
$x^*_1=\ldots=x^*_m=0$. Since all but one (say for $i=1$) of the
sets $\O_i$ are SNC at $\ox$, we have that $\|x^*_{ik}\|\to 0$ as
$k\to\infty$ for $i=2,\ldots,m$. Then \eqref{snc1} implies that
$\|x^*_{1k}\|\to 0$ as well, which verifies implication
(\ref{eq:AQC}) and thus completes the proof of the proposition.
$\h$\vspace*{0.05in}

The following example illustrates the validity of the AQC for
infinite systems of sets.

\begin{Example}[AQC for infinite systems]\label{Ex:AQC} We verify
that the AQC holds in the framework of Example~\ref{Ex:nequi} at
the origin $\ox=(0,0)\in\R^2$. Recall that for each $k\in\N$ the
normal cone to a convex set $\O_k$ from \eqref{ok} at a boundary
point $x=(x_1,x_2)$ is computed by
\begin{equation*}
N(x;\O_k)=\R_+\xi_k(x)\;\mbox{ with }\;\xi_k(x)=\xi_k(x_1,x_2)=\begin{cases}
(2k^m x_1,-1)&\mbox{for }x_1>0,\\
(0,-1)&\mbox{for }x_1\le 0.
\end{cases}
\end{equation*}
If according to the left-hand side of \eqref{eq:AQC} we have
\begin{equation*}
\Big\|\sum_{k\in I_\ve}\lm_{\ve k}\xi_k(x_{\ve k})\Big\|\to
0\;\mbox{ as }\;\ve\dn 0,
\end{equation*}
then it follows from the above representation of $\xi_k$ that its
component goes to zero as $k\to\infty$. Thus
\begin{equation*}
\sum_{k\in I_\ve}\|\lm_{\ve k}\xi_k(x_{\ve k})\|^2\to 0\;\mbox{ as
}\ \ve\dn 0,
\end{equation*}
which verifies the AQC property of the system $\{\O_k\}_{k\in\N}$
at $\ox$.
\end{Example}
\vspace*{0.05in}

Now we are ready to define limiting $R$-normals and derive
infinite intersection rules for them. In the definition below
$R_k$ stands for a rate function for each $x^*_k$; these functions
may be different from each other.

\begin{Definition}[Limiting $R$-normals to infinite set intersections]
Consider an arbitrary set system $\{\O_i\}_{i\in T}\subset X$, and
let $\O:=\bigcap_{i\in T}\O_i$ with $\ox\in\O$. We say that a dual
element $x^*$ is a {\sc limiting $R$-normal} to $\O$ at $\ox$ if
there exist sequences $\{(x_k,x^*_k)\}_{k\in\N}\subset X\times
X^*$ such that $x_k\st{\O}\to\ox$, $x^*_k\st{w^*}{\lto}x^*$ as
$k\to\infty$ and that each element $x^*_k$ is an $R_k$-normal to
$\O$ at $x_k$,
\end{Definition}

It is clear from the definition and Proposition~\ref{Prop:R-Fnor}
that any limiting $R$-normal is a basic/limiting normal to $\O$ at
$\ox$. Conversely, if $T$ is a finite index set and $X$ is an
Asplund space, then we the reverse implication holds, i.e., any
limiting/basic normal is a limiting $R$-normal.\vspace*{0.05in}

The next theorem provides a representation of limiting $R$-normals
to infinite set intersections via Fr\'echet normals to each set
under consideration. In particular, it implies a useful calculus
rule for the basic normal cone \eqref{eq:Lcone} to infinite
intersections.

\begin{Theorem}[Representation of limiting $R$-normals to infinite
intersections]\label{Thm:RepLR} Let $\O:=\bigcap_{i\in T}\O_i$
with $\ox\in\O$ for the system $\{\O_i\}_{i\in T}\subset X$
satisfying the AQC property from Definition~{\rm\ref{Def:AQC}} at
$\ox$. Then for any given limiting $R$-normal to $\O$ at $\ox$ and
any $\ve>0$ we have the inclusion
\begin{equation*}
x^*\in\cl^*\Big\{\sum_{i\in I}x^*_i+\ve\B^*\ \Big|\;x^*_i\in\Hat
N(x_i;\O_i),\;\|x_i-\ox\|<\ve,\;I\subset T\Big\},
\end{equation*}
where $I\subset T$ is a finite index subset. In particular, if all
the limiting/basic normals to $\O$ at $\ox$ are limiting
$R$-normals in this setting, then
\begin{equation}\label{bnc-inf}
N(\ox;\O)\subset\bigcap_{\ve>0}\cl^*\Big\{\sum_{i\in
I}x^*_i+\ve\B^*\ \Big|\; x^*_i\in\Hat N(x_i;\O_i),\;
\|x_i-\ox\|<\ve,\;I\subset T\Big\}.
\end{equation}
\end{Theorem}
{\bf Proof.} Take a sequence $\{x^*_k\}$ of $R$-normals to $\O$ at
$x_k$ with $x_k\to\ox$ and $x^*_k\st{w^*}{\to}x^*$ as
$k\to\infty$. The latter convergence ensures by the Uniform
Boundedness Principle that the set $\{\|x^*_k\|\}_{k\in\N}$ is
bounded in $X^*$. Picking $\ve>0$ sufficiently small, we find
$x_k\in\O$ with $\|x_k-\ox\|<\ve$. Applying Theorem~\ref{Thm:Fuz}
to $x^*_k$ for each $k\in\N$ gives us sequences $x^*_{ik}\in\Hat
N(x_{ik};\O_i)$ with $\|x_{ik}-x_k\|<\ve$ for $i\in I_k\subset T$
and $\lm_k\ge 0$ satisfying
\begin{equation}\label{eq:FuzLR}
\lm_k x^*_k\in\sum_{i\in I_k}x^*_{ik}+\ve\B^*\;\mbox{ and }\;
\lm_k^2+\lm_k^2\|x^*_k\|^2+\sum_{i\in I_k}\|x^*_{ik}\|^2=1,\quad
k\in\N.
\end{equation}
Let us show that the sequence $\{\lm_k\}$ is bounded away from 0.
Assuming on the contrary $\lm_k\dn0$ as $k\to\infty$, we have
\begin{equation*}
\Big\|\sum_{i\in I_k}x^*_{ik}\Big\|\lto 0\;\mbox{as }\ k\to\infty
\end{equation*}
from the inclusion in \eqref{eq:FuzLR}. Then the imposed AQC leads
us to
\begin{equation*}
\sum_{i\in I_k}\|x^*_{ik}\|^2\to 0\;\mbox{ as }\;k\to\infty,
\end{equation*}
which contradicts the equality in (\ref{eq:FuzLR}) and thus shows
that there is constant $C>0$ with $\lm_k>C$ for all $k\in\N$
sufficiently large. Rescaling finally the inclusion in
(\ref{eq:FuzLR}), we get
\begin{equation*}
x^*_k\in\sum_{i\in
I}\frac{x^*_{ik}}{\lm_k}+\frac{\ve}{C}\B^*,\quad k\in\N,
\end{equation*}
which ensures that $x^*_k\st{w^*}{\lto}x^*$ as $k\to\infty$ and
thus justifies the first conclusion of the theorem. The second
ones on basic normals follows immediately. $\h$\vspace*{0.05in}

The next corollary provides more explicit results for the case of
infinite systems of cones, with the replacement of Fr\'echet
normals in Theorem~\ref{Thm:RepLR} by basic normals at the origin.

\begin{Corollary}[Limiting $R$-normals to intersection of cones]
\label{Thm:RepCones} Let $\{\Lm_i\}_{i\in T}$ be a system of cones
in $X$, and let $\Lm:=\bigcap_{i\in T}\Lm_i$. Suppose that $x^*\in
X^*$ is a limiting $R$-normal to $\Lm$ at the origin and that the
AQC property from Definition~{\rm\ref{Def:AQC}} holds at $\ox=0$.
Then for any $\ve>0$ we have the representation
\begin{equation*}
x^*\in\cl^*\Big\{\sum_{i\in I}x^*_i+\ve\B^*\ \Big|\; x^*_i\in
N(0;\Lm_i),\;I\subset T\Big\}
\end{equation*}
via finite index subsets $I\subset T$. If furthermore all the
limiting/basic normals to $\Lm$ at the original are limiting
$R$-normals in this setting, then
\begin{equation*}
N(0;\Lm)\subset\bigcap_{\ve>0}\cl^*\Big\{\sum_{i\in
I}x^*_i+\ve\B^*\ \Big|\; x^*_i\in N(0;\Lm_i),\;I\subset T\Big\}.
\end{equation*}
\end{Corollary}
{\bf Proof.} It is not hard to check that $\Hat
N(w_i;\Lm_i)\subset N(0;\Lm_i)$ for any cone $\Lm_i$ and any
$w_i\in\Lm_i$; see, e.g., \cite[Proposition~2.1]{MorPh10a}. Then
we have both conclusions of the corollary from
Theorem~\ref{Thm:RepLR}. $\h$\vspace*{0.05in}

\begin{Remark}[Comparison with known results]\label{equi} For the
case of finite set systems the intersection rules of
Theorems~\ref{Thm:Fuz} and \ref{Thm:RepLR} go back to the
well-known results of \cite{m-book1}. In fact, not much has been
known for representations of generalized normals to infinite
intersections. Our previous results in this direction obtained in
\cite{MorPh10a,MorPh10b}, obtained on the base of the tangential
extremal principle in finite dimensions, have a different nature
and do not generally reduce to those in \cite{m-book1} for finite
set systems.

An interesting representation of the basic normal cone
\eqref{eq:Lcone} has been recently established in
\cite[Theorem~3.1]{Seidman10} for infinite intersections of sets
given by inequality constraints with smooth functions. This result
essentially exploits specific features of the sets and functions
under consideration and imposes certain assumptions, which are not
required by our Theorem~\ref{Thm:RepLR}. In particular,
\cite[Theorem~3.1]{Seidman10} requires the equicontinuity of the
constraint functions involved, which is not the case of our
Theorem~\ref{Thm:RepLR} as shown in Examples~\ref{Ex:Rnor-inf} and
\ref{Ex:nequi}. Note to this end that all the limiting normals are
limiting $R$-normals in the framework of Example~\ref{Ex:Rnor-inf}
and that the AQC assumption is satisfied therein; see
Example~\ref{Ex:AQC}.
\end{Remark}

We finish the paper with deriving necessary optimality conditions
for problems of semi-infinite and infinite programming with
geometric constraints given by
\begin{equation}\label{sip}
\mbox{minimize}\quad\ph(x)\;\mbox{ subject to}\;x\in\O_i,\quad
t\in T,
\end{equation}
with a general cost function $\ph\colon X\to\oR$ and constraints
sets $\O_t\subset X$ indexed by an arbitrary (possibly infinite)
set $T$. We refer the reader to \cite{CLMP1,GL98,MorPh10b} and the
bibliographies therein for various results, discussions, and
examples concerning optimization problems of type \eqref{sip} and
their specifications. The limiting normal cone representation
\eqref{bnc-inf} for infinite set intersections in
Theorem~\ref{Thm:RepLR}, combined with some basic principles in
constrained optimization, leads us to necessary optimality
conditions for local optimal solutions to \eqref{sip} expressed
via its initial data.

The next theorem contains results of this kind in both {\em lower
subdifferential} and {\em upper subdifferential} forms; see
\cite[Chapter~5]{m-book2} for general frameworks of constrained
optimization and \cite{CLMP1} for semi-infinite/infinite programs
with linear inequality constraints in \eqref{sip}. The lower
subdifferential condition is given below for the case of locally
Lipschitzian cost functions on Asplund spaces via the construction
\begin{equation*}
\partial\ph(\ox):=\Limsup_{x\to\ox}\Hat\partial\ph(x)
\end{equation*}
known as the {\em Mordukhovich/basic/limiting subdifferential} of
$\ph$ at $\ox$; see
\cite{Borwein-Zhu-TVA,m-book1,Rockafellar-Wets-VA,sc} for more
details and discussions. The upper subdifferential condition below
employs the so-called {\em Fr\'echet upper
subdifferential/superdifferential} of $\ph$ at this point defined
by
\begin{equation*}
\Hat\partial^+\ph(\ox):=-\Hat\partial(-\ph)(\ox).
\end{equation*}

\begin{Theorem}[Necessary optimality condition for semi-infinite and infinite
programs with general geometric constraints] Let $\ox$ be a local
optimal solution to problem \eqref{sip}. Assume that any basic
normal to $\O:=\bigcap_{i\in T}\O_i$ at $\ox$ is a limiting
$R$-normal in this setting, and that the AQC requirements is
satisfied for $\{\O_i\}_{i\in T}$ at $\ox$. Then the following
conditions, involving finite index subsets $I\subset T$, hold:

{\bf (i)} For general cost functions $\ph$ finite at $\ox$ we have
\begin{equation}\label{up}
-\Hat\partial\ph(\ox)\subset\bigcap_{\ve>0}\cl^*\Big\{\sum_{i\in
I}x^*_i+\ve\B^*\ \Big|\ x^*_i\in\Hat N(x_i;\O_i),\;
\|x_i-\ox\|<\ve,\;I\subset T\Big\}.
\end{equation}

{\bf (ii)} If in addition $\ph$ is locally Lipschitzian around
$\ox$, then
\begin{equation}\label{lo}
0\in\partial\ph(\ox)+\bigcap_{\ve>0}\cl^*\Big\{\sum_{i\in
I}x^*_i+\ve\B^*\ \Big|\ x^*_i\in\Hat N(x_i;\O_i),\;
\|x_i-\ox\|<\ve,\;I\subset T\Big\}.
\end{equation}
\end{Theorem}
{\bf Proof.} It follows from \cite[Proposition~5.2]{m-book2} that
\begin{equation}\label{u-opt}
-\Hat\partial\ph(\ox)\subset\Hat N(\ox;\O)\subset N(\ox;\O)
\end{equation}
for the general constrained optimization problem
\begin{equation}\label{opt}
\mbox{minimize }\;\ph(x)\;\mbox{ subject to }\;x\in\O.
\end{equation}
Employing now in \eqref{u-opt} the intersection formula
\eqref{bnc-inf} for basic normals to $\O=\bigcap_{i\in T}\O_i$, we
arrive at the upper subdifferential necessary optimality condition
\eqref{up} for problem \eqref{sip}.

To justify \eqref{lo}, we get from \cite[Propostion 5.3]{m-book2}
the lower subdifferential necessary optimality condition
\begin{equation}\label{lo-opt}
0\in\partial\ph(\ox)+N(\ox;\O)
\end{equation}
for problem \eqref{opt} provided that $\ph$ is locally
Lipschitzian around $\ox$. Using the intersection formula
\eqref{bnc-inf} in \eqref{lo-opt} completes the proof of the
theorem. $\h$

\end{document}